\documentclass[a4paper,11pt]{article}
\pdfoutput=1 
\usepackage{setspace}

\newif\ifCOMMENTS
\COMMENTSfalse

\usepackage{amssymb}

\usepackage{natbib}

\usepackage{fullpage} 
\usepackage{amsmath} 
\usepackage{amsmath,calc} 
\usepackage{adjustbox}
\usepackage[normalem]{ulem}
\usepackage{color}
\newlength{\LPlhbox}

\usepackage{algorithm} 
\usepackage{algorithmic} 
\renewcommand{\algorithmiccomment}[1]{\bgroup\hfill//~#1\egroup}  

\usepackage{float} 

\usepackage{bbm} 
\usepackage{comment}
\usepackage{tabularx} 
\newcolumntype{Y}{>{\centering\arraybackslash}X} 
\usepackage{multirow}
\usepackage{subcaption}
\usepackage{paralist}

\allowdisplaybreaks
\usepackage{graphicx}
\usepackage{lscape}

\usepackage[dvipsnames]{xcolor}
\usepackage{diagbox} 
\usepackage{slashbox}

\usepackage[colorlinks=false, urlcolor=blue, linkcolor=false]{hyperref} 
\usepackage{url}

\usepackage{pifont}
\usepackage{booktabs}
\usepackage[capitalize,noabbrev]{cleveref}
\usepackage{xspace}
\newcommand{\probabbr}{HHCRSP-TS\xspace}
\newcommand{\probabbrtt}{HHCRSP-TS-TT\xspace}
\newcommand{\Tmax}{T_{\mathrm{max}}}

\newcommand{\argmin}{\mathrm{arg\,min}}
\newcommand{\Astart}{A_{\mathrm{start}}}
\newcommand{\Acare}{A_{\mathrm{care}}}

\newcommand{\LGraph}{\mathcal{G}}
\newcommand{\LNodes}{\mathcal{N}}
\newcommand{\LNodesVisit}{\LNodesV{V}}
\newcommand{\LNodesV}[1]{\LNodes(#1)}

\newcommand{\LArcs}{\mathcal{A}}
\newcommand{\LAstart}{\LArcs_{\mathrm{start}}}
\newcommand{\LAcare}{\LArcs_{\mathrm{care}}}
\newcommand{\LAwait}{\LArcs_{\mathrm{wait}}}

\newcommand{\TINOSPLIT}{\texttt{TI$^-$}\xspace}
\newcommand{\TIALLSPLIT}{\texttt{TI$^+$}\xspace}
\newcommand{\TI}{\texttt{TI}\xspace}
\newcommand{\TIHMTZ}{\texttt{TI+HMTZ}\xspace}
\newcommand{\TIHTI}{\texttt{TI+HTI}\xspace}

\newcommand{\tu}{t_{u}}
\newcommand{\tv}{t_{v}}
\newcommand{\tw}{t_{w}}

\usepackage{enumerate}

\defcitealias{GVR2023}{Grand View Research (GVR), 2023}
\defcitealias{COAVVT2023}{SOVVT (2023)}
\usepackage{appendix}

\author{Loek van Montfort\thanks{l.g.a.j.van.montfort2@vu.nl}, Wout Dullaert\thanks{w.e.h.dullaert@vu.nl}, and Markus Leitner\thanks{m.leitner@vu.nl}}
\date{\textit{Department of Operations Analytics, Vrije Universiteit Amsterdam, The Netherlands}\\[2ex]%
    May 28, 2025}

\title{Task-splitting in home healthcare routing and scheduling}

\begin{document}
\maketitle
\begin{abstract}
This paper introduces the concept of task-splitting into home healthcare (HHC) routing and scheduling. It focuses on the design of routes and timetables for caregivers providing services at patients' homes. Task-splitting is the division of a (lengthy) patient visit into separate visits that can be performed by different caregivers at different times. The resulting split parts may have reduced caregiver qualification requirements, relaxed visiting time windows, or a shorter/longer combined duration. However, additional temporal dependencies can arise between them.
To incorporate task-splitting decisions into the planning process, we introduce two different mixed integer linear programming formulations, a Miller–Tucker–Zemlin and a time-indexed variant. These formulations aim to minimize operational costs while simultaneously deciding which visits to split and imposing a potentially wide range of temporal dependencies. We also propose pre-processing routines for the time-indexed formulation and two heuristic procedures. These methods are embedded into the branch-and-bound approach as primal and improvement heuristics. 
The results of our computational study demonstrate the additional computational difficulty introduced by task-splitting possibilities and the associated additional synchronization, and the usefulness of the proposed heuristic procedures. 
From a planning perspective, our results indicate that integrating task-splitting decisions into the planning process reduces staff requirements, decreases HHC operational costs, and allows caregivers to spend relatively more time on tasks aligned with their qualifications. 


\end{abstract}
\section{Introduction}\label{s: introduction}
In industrialized countries, healthcare providers are confronted with the consequences of an aging population, which drives an increase in the demand for healthcare services. One of the sectors particularly affected by this is home healthcare (HHC), which will play an increasingly prominent role in making healthcare affordable. In the United Kingdom, for instance, this translates into a projected compound annual growth rate for the sector of 9.1\% until 2030 \citepalias{GVR2023}. HHC involves caregivers traveling to patients' homes to provide various services, including assistance with taking medication, administering injections, wound care, or showering. The benefits of HHC are twofold. Firstly, people prefer to age at home \citep{Wiles2011} as it allows them to retain their autonomy. 
Secondly, it is a cost-effective alternative to traditional hospital care. However, the supply of HHC services is already under pressure due to staff shortages, so meeting the growing demand from an aging population presents a significant challenge. 

Driven by the social relevance of affordable and reliable HHC, home healthcare planning has received significant research attention in recent years. The focus has been on the design of efficient routes and corresponding schedules for caregivers, formally known as Home Healthcare Routing and Scheduling Problems (HHCRSPs). 
These planning problems typically include issues such as: i) caregivers having specific capacities and work preferences, ii) patient visits having service and time constraints, and iii) the need to synchronize multiple visits for the same patient. Given the complexity of these planning processes, further optimization of caregiver schedules is one element that HHC providers can do to improve caregiver allocation. 

\medskip 

The diverse national and regulatory settings have led to the emergence of numerous variants of the HHCRSP. Overviews of the wide range of modeling and solution approaches utilized are provided by \citet{Masmoudi2024,EUCHI2022, DIMASCOLO2021107255, FIKAR201786}. Most research relates to the daily planning of care, with time windows indicating when care can be provided. Other commonly considered characteristics are the required qualifications of visiting caregivers, the permitted working hours, and the legislative rules for caregivers. 

A significant part of the literature related to HHCRSPs does not consider the temporal restrictions 
between patient visits (see, e.g., \cite{PAHLEVANI2022102878, LI2021102420, BRAEKERS2016428}) that further complicate the already challenging task of solving the basic types of HHCRSPs. Temporal dependencies create interdependence between caregiver routes and schedules, resulting in complex planning challenges \citep{SOARES2024817, Drexl2012}. For example, strict synchronization occurs if two caregivers are required for a certain task, such as lifting a disabled person. Precedence refers to a certain minimum and/or maximum amount of time needed between two visits (e.g., medication before or after dinner). Disjunctive requirements imply that two tasks must be executed on non-overlapping moments in time (e.g., assistance with bathing and a medical treatment).
These temporal dependencies are more generally referred to as synchronization and may also arise when a patient's care tasks are split among multiple visits. 
Synchronization restrictions can be required for approximately 20\% of the patient visits \citep{Polnik2020}, triggering a growing academic interest on how to handle them. These research efforts are supported by the development of a wide variety of heuristic approaches, see, e.g., \citep{Masmoudi2023, Bredstrom2008} for the articles considering the simultaneous execution of tasks and \citep{OLADZADABBASABADY2023105829, Frifita2020, Mankowska2013} for works considering other temporal restrictions. 

On the contrary, only a few exact solution methods for HHCRSPs with synchronization requirements have been proposed. 
\cite{RASMUSSEN2012698} developed a Branch-and-Price (BP) approach for an HHCRSP with strict synchronization and precedence. Their BP utilizes a set partition model and enforces synchronization via time window branching. 
\cite{Mankowska2016} proposes a combinatorial Benders approach for solving an HHCRSP with strict synchronization and precedence. Their master problem fixes the caregiver routes, and the continuous subproblem determines the timing of the route (e.g., temporal dependencies). \cite{Qiu2021} describe a Branch-and-Price-and-Cut (BPC) for an HHCRSP where the starting times of synchronized visits must be within a certain range of each other. 
They create initial caregiver routes with Adaptive Large Neighborhood Search and solve the pricing problem with a specially designed labeling algorithm. 
\cite{DOULABI2020} investigate an HHCRSP that includes strict synchronization and stochastic travel and service times. They solve it with an L-shaped and a Branch-and-Cut algorithm, using valid inequalities for subtours, capacity constraints, and no overlap constraints. Despite the ability of exact solution methods to provide optimal solutions, their practical use is limited by high computational needs. 

In addition to the research on HHCRSPs were caregivers travel independently, settings are studied in which teams of two caregivers travel together. Synchronization can then be used to allow the splitting of some caregiver teams traveling together. In \cite{DEAGUIAR2023101503, Nozir2020}, patient visits requiring a team of two caregivers are allowed to be also served by two simultaneous visits of caregivers traveling independently. The results indicate that this can result in less staff being required and higher planning efficiency, suggesting the potential of splitting resources.

Studies on split-delivery Vehicle Routing Problems (VRPs) allow for the splitting of tasks or demand to enhance planning flexibility. Contrary to health-care settings, split-delivery VRPs do, however, typically consider arbitrary splits of the total amount of individual customers. One notable extension of this concept of arbitrary splitting has been proposed recently by \cite{FERNANDEZ20181078} who study a multi-depot VRP with pre-defined splitting options. 

A broader perspective reveals that many HHCRSPs can be classified as specific types of VRPs with synchronization. Such problems have been addressed by \citet{SHI2020124112, LIU2019250, Parragh2018, Afifi2016}. While the majority of research on exact solution methods for VRPs with synchronization focuses on a homogeneous fleet \citep{DOHN2011, BREDSTROM2007}, there are also articles considering heterogeneous sets of vehicles. For instance, \citet{Luo2016,DOHN20091145} developed BP(C) algorithms for set partition formulations, securing synchronization through time window branching.

In practice, HHC providers have already explored various methods for organizing their work more efficiently. In the Netherlands, for example, some organizations use self-managing teams in which each caregiver typically performs all different care tasks required during a single patient visit (Buurtzorg). 
Another approach consists of differentiating the care tasks requested by a patient and assigning them to caregivers with appropriate skill levels. 
This approach enables caregivers with specific skills to perform few tasks outside their area of expertise, which can help to reduce specific staff shortages.

Although assigning care tasks to patient visits is part of the HHC planning process, the current literature on HHCRSP considers the visits as given, implicitly assuming that these assignment decisions are already made in advance. As a consequence, the question whether and when to divide care tasks and how to assign them to different types of caregivers has not been formally investigated yet. 
Furthermore, no tools exist that support such decisions regarding the organization of care task into visits while taking into account caregiver schedules and routes. 
Integrating these decisions, instead of deciding which tasks to divide in advance, can support the identification of visits for which splitting tasks is beneficial for the overall planning. 
Motivated by these observations and discussions with practitioners from the healthcare sector, this paper formally investigates whether and when \emph{task-splitting} can contribute to enhanced caregiver planning in the home healthcare sector.

\medskip 

Task-splitting is the division of a (potentially inefficient) lengthy patient visit into multiple, shorter visits. This can introduce additional scheduling flexibility for HHC providers and eventually help to obtain better plannings. To the best of our knowledge, this is the first paper to incorporate the decision of whether a visit should be split into two parts within the planning process. Combining task-splitting with routing and scheduling decisions ensures that the proposed solution algorithm supports decision making about the use of task-splitting while considering the eventual effect on HHC schedules and routes. 
The advantages of task-splitting emerge in a number of ways. First, it can result in a \emph{better use of the available staff}, as one of the split parts may require lower caregiver qualifications. Second, it allows the \emph{split parts to be performed at different times}, thus reducing the peak demand. Third, it can result in \emph{improved alignment between patient needs and caregivers}, as caregiver qualifications can be indicated for both split parts. Lastly, task-splitting can lead to \emph{efficiency gains}, for example, by eliminating mandatory waits between two tasks. 
Caution should, however, be taken when deciding whether to split tasks. For instance, patients may not wish to be visited from multiple caregivers or be forced to remain at home all day. Moreover, splitting tasks can introduce additional working time for caregivers (e.g., additional travel or start-up time). Temporal restrictions between the execution of the different split parts must also be taken into account, which complicates the planning problem. For example, tasks requiring contact with the patient, cannot be performed simultaneously. This type of non-overlapping temporal restriction is another element that has not been addressed in the current exact solution methods for HHCRSPs. 
The potential of task-splitting extends beyond the scope of HHC. It could also occur in care homes or hospitals, for example. Therefore, the relevance of this research is broader than home care; it is also of interest to related (healthcare) settings that face similar (future) staffing challenges. 

\medskip 
This paper aims to support the incorporation of task-splitting into home healthcare planning and to explore its impact on HHC providers and related settings. Its main contributions are:

\begin{itemize}
    \item We propose the first variant of the HHCRSP that incorporates task-splitting possibilities. While our problem variant considers a wide range of temporal dependency constraints, caregiver qualifications, and operational costs, it does not incorporate other, practically relevant aspects (such as, e.g., client preferences, overtime or soft time windows) that could be integrated in future work. By focusing on synchronization and the possibility of splitting visits into two separate parts within a fundamental HHCRSP, we can develop an exact solution approach that allows us to assess the impact of task-splitting in HHC routing and scheduling. 
    \item We develop two integer linear programming formulations incorporating task-splitting decisions: a Miller-Tucker-Zemlin formulation and a time-indexed formulation.
    \item We introduce pre-processing techniques and two heuristic procedures to improve the computational performance of the solution algorithm. The algorithm supports decision making about the use of task-splitting considering the eventual effect on HHC schedules and routes.
    \item We explore the computational challenges and managerial benefits of integrating the possibility of task-splitting in HHC planning in various scenarios. 
\end{itemize}

The remainder of this paper is organized as follows. Section \ref{s: problem_description} defines the new problem variant. Section \ref{sec: math_form} introduces a Miller-Tucker-Zemlin and a time-indexed integer linear program (ILP). Section \ref{s: sol_approach} discusses pre-processing techniques and proposes two heuristic procedures to enrich an ILP solver. Section \ref{sec:instances} describes the generation of benchmark instances. \cref{sec:results} presents the results of our computational study, which incorporates an analysis of the computational characteristics and an exploration of the effect of task-splitting on caregiver allocations in various scenarios. Finally, we summarize the findings and outline potential directions for further research in Section \ref{sec:conclusion}. 

\section{Problem definition}\label{s: problem_description}
This section formally introduces the Home Healthcare Routing and Scheduling Problem with Task-Splitting (\probabbr), which aims to design optimal routes and timetables for caregivers that provide a service at patients' homes. 
The problem considers a set of available \emph{caregivers}, denoted by $C$, who have to perform a set of \emph{original patient visits} $V^{\mathrm{o}}=\{1, \dots, n\}$ within a planning horizon $[0, T_{\mathrm{max}}]$. Each caregiver $k\in C$ is associated with a qualification type $q_k$ from the set of qualification types $Q$ that represent, e.g., qualification levels or certain specializations and is available throughout the whole planning horizon. Each patient visit $v\in V^{\mathrm{o}}$ is associated with a set of qualifications $Q_v\subseteq Q$ that specify which types of caregivers can perform visit $v$, the visit duration $d_v$, $0<d_v\le T_\mathrm{max}$, and a time window $[\alpha_{v},\beta_{v}]$, $0\le \alpha_v < \beta_v \le T_{\mathrm{max}}-d_v$, during which visit $v$ must start. A patient can require multiple visits.

The subset $S \subseteq V^{\mathrm{o}}$ contains \emph{visits} that are comprised of multiple tasks that can be split into two separate parts. Each \emph{splittable visit} $v \in S$ can either be performed as a whole or replaced by two split visits $v_p$, $p\in \{1,2\}$ with qualification types $Q_{v_p}\subseteq Q$, durations $d_{v_p}$, $0<d_{v_p}\le T_\mathrm{max}$, and time windows $[\alpha_{v_p},\beta_{v_p}]$, $0\le \alpha_{v_p} < \beta_{v_p} \le T_{\mathrm{max}}-d_{v_p}$.
When combined with the unsplit visits, this generates the set of potential patient visits $V = V^{\mathrm{o}} \cup \{v_{p} : v \in S, \; p \in \{1,2\}\}$. 
When the splittable visit $v \in S$ is split into two separate parts, one of the split parts $v_p$, $p\in \{1,2\}$ may be performed by a wider range of caregivers (i.e., $ Q_{v} \subseteq Q_{v_{p}}$). Moreover, the combined duration of the split parts may differ from that of the original visit $v$ due to efficiency gains or losses (i.e., $d_{v} \neq d_{v_1} + d_{v_2}$). It may also be possible to perform the split parts at different points in time (i.e., $[\alpha_{v_{p}},\beta_{v_{p}}] \neq  [\alpha_v,\beta_v]$).

Set $D\subseteq V^2$ encompasses all pairs of visits with temporal dependencies, such as precedence requirements, overlapping or non-overlapping visiting times, or a simultaneous start of the visits. The \probabbr considers a generic set of \emph{synchronization requirements} by associating parameters $\delta_{uv}^{\mathrm{min}}$ and $\delta_{uv}^{\mathrm{max}}$ to each ordered pair $(u,v)$ of visits $\{u,v\}\in D$ with temporal dependencies. This yields four parameters for each such pair of visits. Thus, $\delta_{uv}^{\mathrm{min}}$ and $\delta_{uv}^{\mathrm{max}}$ indicate the permitted minimum and maximum difference in the starting times of visits $u$ and $v$ if visit $u$ starts before or simultaneously with visit $v$. This means that starting time $t \in [\alpha_{u}, \beta_{u}]$ imposes the start of visit $v$ within the time window $[t+\delta_{uv}^{\mathrm{min}},t+\delta_{uv}^{\mathrm{max}}]$. \cref{tab:synchronization:overview} summarizes various synchronization types and their corresponding parameter values. 

\renewcommand{\arraystretch}{0.65}
\begin{table}
	\caption{Considered synchronization types and corresponding parameter values. Here, $u \preceq v$ indicates that visit $u$ must start before visit $v$, while $\Delta_{\mathrm{min}}$ and $\Delta_{\mathrm{max}}$ specify the minimal and maximal starting time differences.}
	\label{tab:synchronization:overview}
 \centering 
 \small
\begin{tabular}{cc cc cc}
	\toprule
	\multicolumn{2}{c}{synchronization type} & \multicolumn{4}{c}{parameter values} \\
	\cmidrule(lr){1-2} \cmidrule(lr){3-6}	
	ordering & temporal aspect & $\delta_{uv}^{\mathrm{min}}$ & $\delta_{uv}^{\mathrm{max}}$ & $\delta_{vu}^{\mathrm{min}}$ & $\delta_{vu}^{\mathrm{max}}$  \\
	\midrule
	- & strict synchronization (simultaneous start) & 0 & 0 & 0 & 0 \\	
	- & limit difference of start (min.\ $\Delta_{\mathrm{min}}$, max.\ $\Delta_{\mathrm{max}}$) & $\Delta_{\mathrm{min}}$ & $\Delta_{\mathrm{max}}$ & $\Delta_{\mathrm{min}}$ & $\Delta_{\mathrm{max}}$ \\
	-& visits must overlap & 0 & $d_u$ & 0 & $d_v$ \\
    -& visits must not overlap & $d_u$ & $\Tmax$ & $d_v$ & $\Tmax$ \\
	\midrule
	$u \preceq v$ & visit must be completed   & $d_u$ & $\Tmax$ & $\Tmax$ & $\Tmax$ \\
	$v \preceq u$ & visit must be completed & $\Tmax$ & $\Tmax$ & $d_v$ & $\Tmax$ \\	
	\midrule
	$u \preceq v$ &  limit difference of start (min.\ $\Delta_{\mathrm{min}}$, max.\ $\Delta_{\mathrm{max}}$) & $\Delta_{\mathrm{min}}$ & $\Delta_{\mathrm{max}}$ & $\Tmax$ & $\Tmax$ \\
	$v \preceq u$ &limit difference of start (min.\ $\Delta_{\mathrm{min}}$, max.\ $\Delta_{\mathrm{max}}$)   & $\Tmax$ & $\Tmax$ & $\Delta_{\mathrm{min}}$ & $\Delta_{\mathrm{max}}$ \\
	\bottomrule
\end{tabular}
\end{table}

Parameter $t_{uv}\ge 0$ indicates the (potentially asymmetric) time it takes a caregiver to travel from visit $u\in V$ to visit $v\in V$  while cost parameter $e^{q} \ge 0$ represents the wage (per time unit) of a caregiver with qualification $q\in Q$. The objective of the \probabbr is to identify a plan with \emph{minimal overall costs} in which each patient visit is performed by a sufficiently qualified caregiver within the given time window and in which all temporal dependencies are respected. For each splittable visit $v\in S$, either the original visit $v$ or both split visits $v_p$, $p\in \{1,2\}$, must be performed. The costs of each plan are equal to the sum of all total working times of all caregivers multiplied by the wage parameter commensurate with each caregiver's qualification. Thereby, the total working time of a caregiver is defined as the sum of time spent on patient visits, the travel time between visits, and any waiting time possibly required before or after visits in order to meet the time windows of subsequent visits. The problem, thus, incorporates for each caregiver a wage from the start of their first visit until the completion of their last visit. Let $\mathcal{P}$ denote the set of all feasible plans in the \probabbr, $t_{\mathrm{start}}^k(P)$ denote the starting time of the first visit and $t_{\mathrm{end}}^k(P)$ donate the end time of the last visit of caregiver $k\in C$ in plan $P \in \mathcal{P}$. Then, the objective of the \probabbr is to identify a plan $P^*$ that minimizes the total cost, i.e.,
\[
 P^*\in \argmin_{P\in \mathcal{P}} \sum_{k\in C} e^{q_k} (t_{\mathrm{end}}^k(P) - t_{\mathrm{start}}^k(P)). \label{eq:objective}
\]

\begin{figure}[h]
\centering 
\includegraphics[width=.8\linewidth]{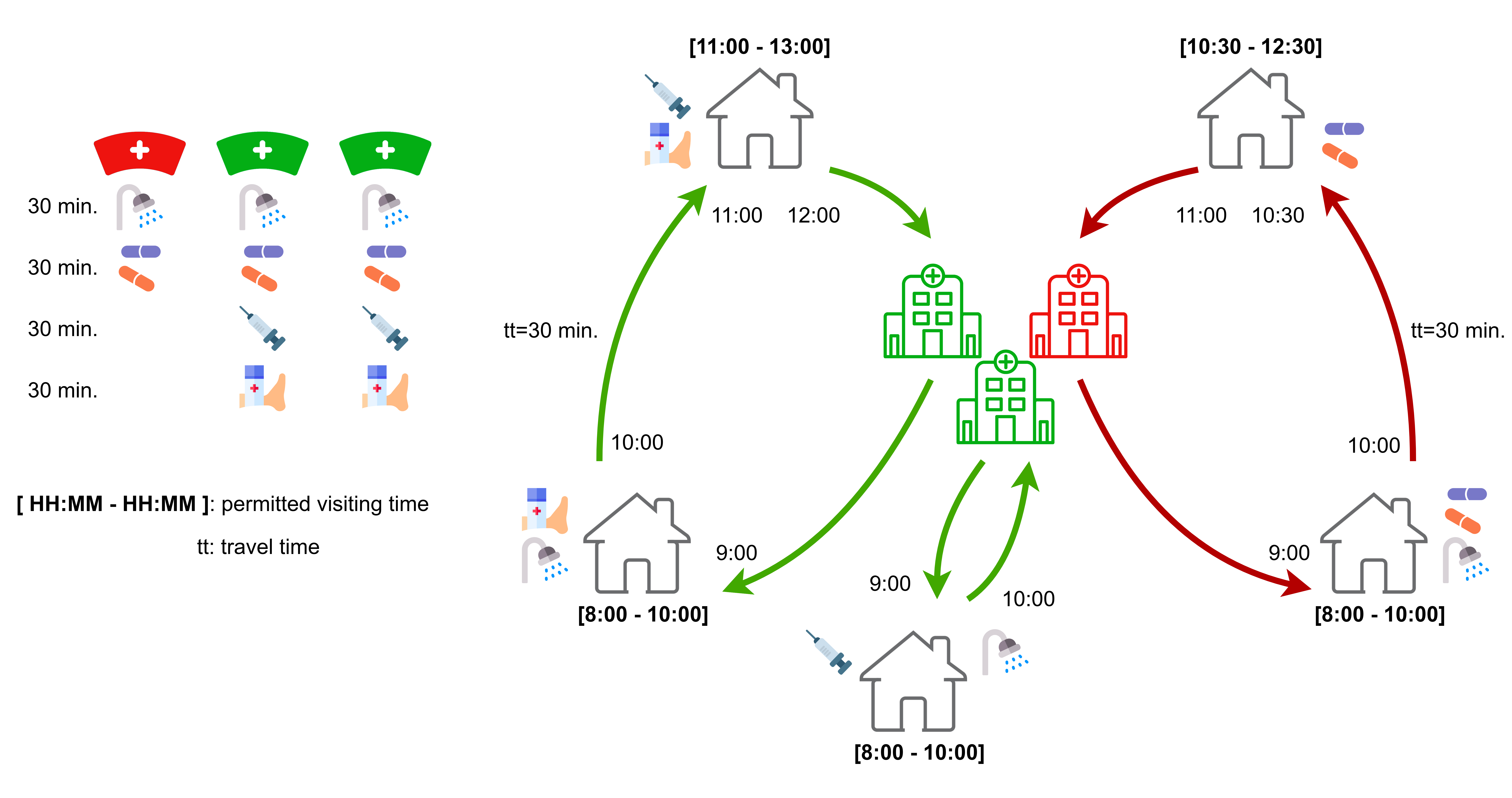}
\caption{Visualization of a \probabbr instance and a corresponding solution without task-splitting.}
\label{fig:example:tasksplitting:nosplit}
\end{figure}

\begin{figure}[h]
\begin{subfigure}[t]{.5\linewidth}
\includegraphics[width=\linewidth]{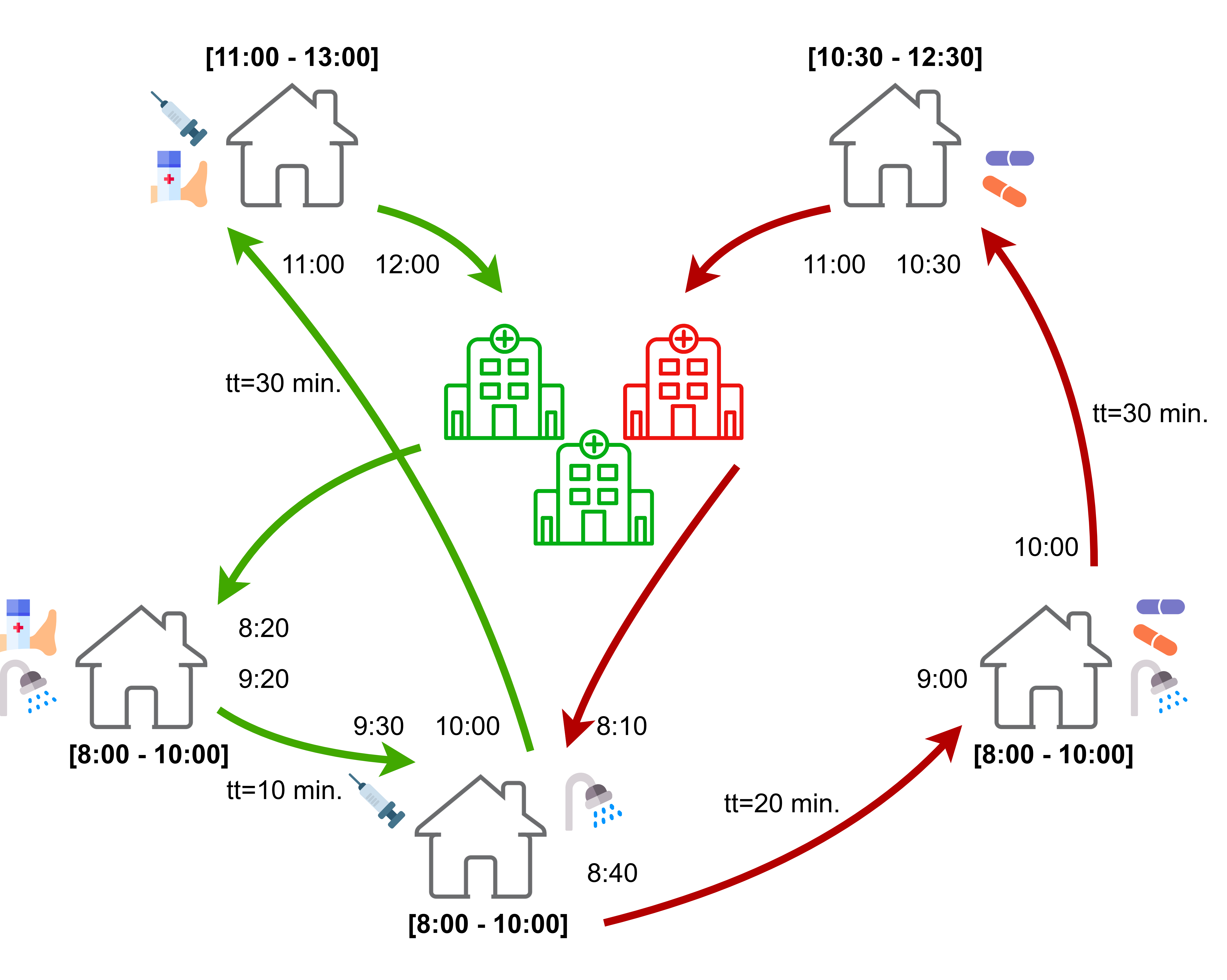}
\caption{fewer caregivers}
\label{fig:example:tasksplitting:split:a}
\end{subfigure}
\begin{subfigure}[t]{.5\linewidth}
\includegraphics[width=\linewidth]{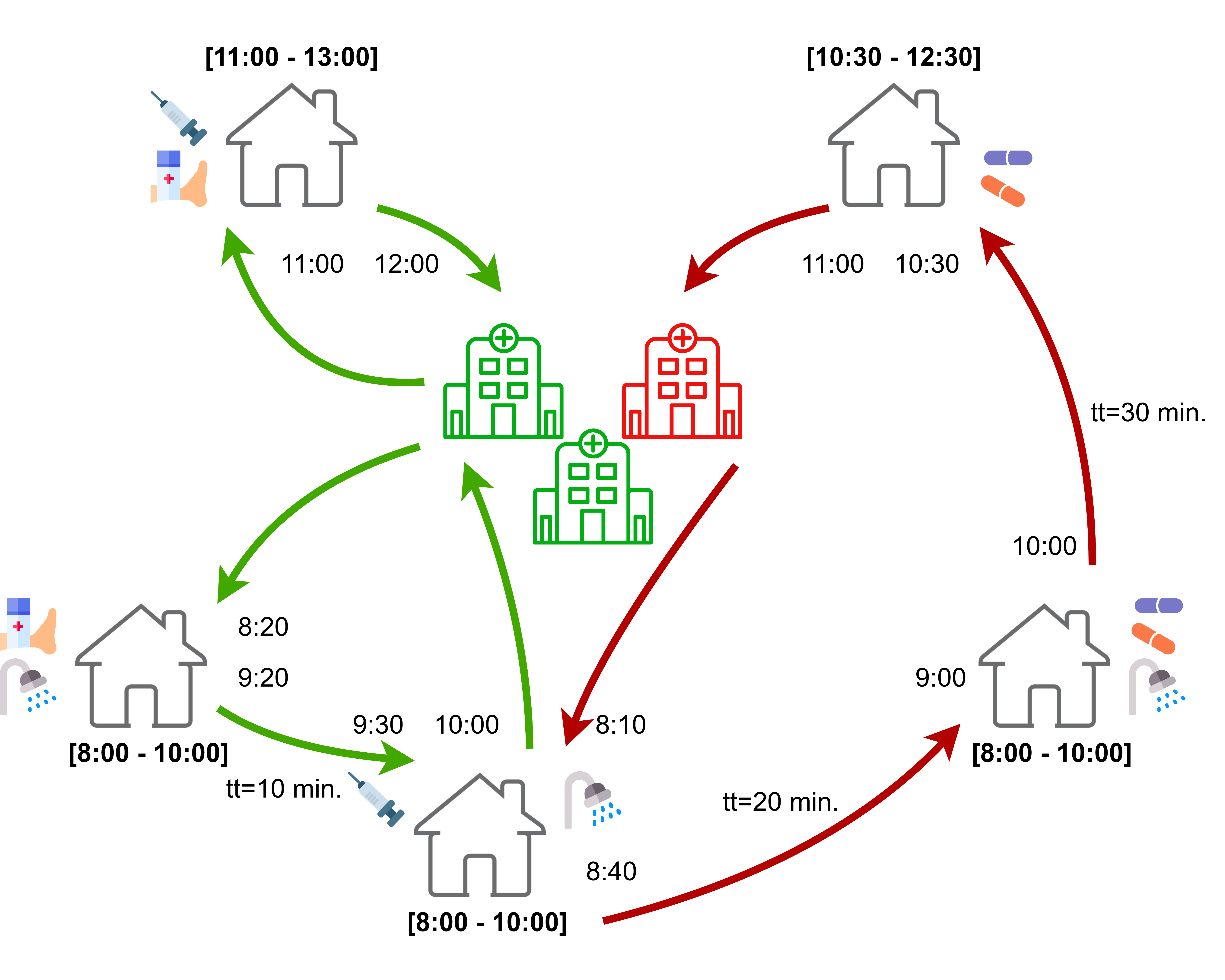}
\caption{less time}
\label{fig:example:tasksplitting:split:b}
\end{subfigure}
\caption{Visualization of two feasible solutions for the \probabbr instance shown in \cref{fig:example:tasksplitting:nosplit}, using task-splitting to provide care with (a) fewer caregivers or (b) in less time.}
\label{fig:example:tasksplitting:split}
\end{figure}

An example of the \probabbr\ that illustrates a situation involving three visits in the early morning and two at the end of the morning is provided in \cref{fig:example:tasksplitting:nosplit}. Each house represents the location of an original patient visit while every hospital indicates the artificial initial (and final) location of a caregiver. The interval between the earliest starting time and latest finishing time of the care of a visit (i.e., the permitted visiting time) is indicated within the associated black brackets and the requested care tasks are indicated by icons. Introducing the option to split the tasks of the visit at the house at the bottom means both tasks no longer have to be performed immediately after each other by the most medically trained type of caregiver. 
As shown in \cref{fig:example:tasksplitting:split}, this allows all care requests to be fulfilled with either fewer caregivers (2 instead of 3) or in less time (-30 minutes). In this example, task-splitting results in shifting work to the least medically trained caregiver. Note that splitting tasks for one of the other visits cannot reduce the number of required caregivers or the total working hours. This indicates a potential drawback of dividing and differentiating care task already in advance.
 
\smallskip

Since related research focuses, for example, on the travel time of caregivers (see, e.g., \cite{PAHLEVANI2022102878, Frifita2020}), we also consider a variant of the \probabbr whose objective is to minimize the total travel time of all caregivers (thus ignoring visit time, waiting times, and wage differences). \ref{app:travel_time_minimization} discusses this problem variant and compares the potential impact of task-splitting for both objective functions.

\paragraph{Assumptions and notation}
In the remainder of this paper, we assume that the same caregiver cannot perform two split parts of a splittable visit consecutively, as this would correspond to an unsplit visit. We also assume that the travel time between patient visits satisfies the triangle inequality and that caregivers always travel to their next patient before waiting. The latter assumption is without loss of generality since the objective is indifferent to the fact whether waiting between two consecutive visits occurs before or after traveling to the next visit. 
To simplify the notation, we add the artificial visits $0$ and $n+1$ to model the start and end of a workday, respectively, resulting in the total set of visits $N=V\cup \{0,n+1\}$. All travel times from or to these artificial nodes are assumed to be equal to zero, i.e., $t_{0v}=0$, $\forall v\in N\setminus \{0\}$ and $t_{v,n+1}=0$, $\forall v\in N\setminus \{n+1\}$. We also assign a time window $[0,\Tmax]$ to both artificial nodes that can be visited by caregivers of all qualifications, i.e., $\alpha_0=\alpha_{n+1}=0$, $\beta_0=\beta_{n+1}=\Tmax$, and $Q_0=Q_{n+1}=Q$. 
Throughout this paper, we use the following notation. The set of all caregivers with qualification $q\in Q$ will be denoted as $C_q=\{k\in C : q_k=q\}$, $C_q\subseteq C$. 
For a directed graph $G'=(V',A')$ and a set of nodes $W\subset V'$, we use $\delta^+(W)=\{(u,v)\in A'\mid u\in W, v\notin W\}$. For simplicity, we write $\delta^+(u)$ instead of $\delta^+(\{u\})$ for node sets of size one. Finally, for a set of variables $Y$ defined on set $U$ and subset $U' \subseteq U$, we use the notation $Y(U')=\sum_{u\in U'} Y_u$.

\section{Mathematical formulations}\label{sec: math_form}  
This section presents two mixed-integer linear programming formulations for the \probabbr. The formulation of \cref{sec:mtz} uses Miller-Tucker-Zemlin (MTZ) constraints \citep{Miller1960}. It is well known, however, that such formulations often fail to solve non-trivial instances of a given problem due to their weak dual bounds. Therefore, \cref{sec:lg} introduces an alternative formulation based on time-indexed flows \citep{Picard1978}. Time-indexed formulations have been successfully used for modeling and solving VRPs with synchronization constraints \citep{FINK2019699,DOHN2011} and a range of other settings, see, e.g., \citep{GOUVEIA2017908,Gouveia2019layered}. Both formulations are used in a combined way in the solution approaches described in \cref{sec:heuristics}.

\subsection{Miller-Tucker-Zemlin formulation}\label{sec:mtz}

The formulation introduced in this section is defined on the directed graph $G=(N,A)$, whose nodes $N$ correspond to all potential visits and whose arcs $A \subset N \times N$ represent potential travels between them. Arc set $A$ is the disjoint union of starting arcs $\Astart=\{(0,v) : v\in V\}$, which indicate travel to a patient's location at the start of a workday, and care arcs $\Acare=\{(u,v) : u\in V, v\in V\cup \{n+1\}, u\ne v, Q_u\cap Q_v\ne \emptyset, \alpha_u + d_u + t_{uv} \le \beta_v\}\setminus \{(u,v) : \exists w \in S : u,v \in \{w,w_1,w_2\}\}$. Each care arc $(u,v)\in \Acare$ represents performing visit $u\in V$ and then traveling to visit $v\in V\cup \{n+1\}$. 
Note that care arcs between two visits are only included if they can be performed by a single caregiver in the given order. The latter is only possible if at least one qualification type allows both to be performed and one can arrive at the second visit before its latest starting time. Note also that travels between a splittable visit and its two split visits are excluded, since either the split visits or the original visit are performed, but not both. Similarly, direct travels between two split visits of the same original visit are excluded because the two split visits may not be performed immediately after each other (in this case, the original visit is performed).

\medskip 

Formulation~\eqref{eq:mtz} uses the following sets of variables. For each arc $a=(u,v)\in A$ and feasible qualification $q\in Q_u\cap Q_v$, travel variables $x_a^q\in \{0,1\}$ indicate whether a caregiver with qualification $q$ travels along arc $(u,v)$, i.e., performs visits $u$ and $v$ immediately after each other. Split variables $s_v\in \{0,1\}$ indicate whether a splittable visit $v\in S$ is split (in which case split parts $v_1$ and $v_2$ must be performed) or not (in which case the original visit $v$ is executed) while variables $y_v\in \{0,1\}$ indicate for each possible type of visit $v \in V$ whether it is performed or not (e.g., a split part). 
Precedence variables $p_{uv}\in \{0,1\}$ are used to model the sequence in which pairs of visits $\{u,v\}\in D$ with temporal dependencies are performed. Thereby, pairs of visits $\{u,v\}\in D$ are considered in an ordered manner such that $u<v$, and $p_{uv}=1$ indicates that visit $u$ starts before or at the same time as visit $v$. Conversely, $p_{uv}=0$ indicates either that visit $v$ starts before or at the same time as visit $u$, or that at least one of these two visits is not performed at all (which can only happen for splittable visits and their split parts). Continuous variables $b_v\ge 0$ represent the starting time of the care of visit $v\in V$ and are set to zero for splittable or split visits that are not performed. Additionally, if a visit $v \in V$ is the first visit of a caregiver route of type $q \in Q_v$, then variable $k_v^q$ is equal to the starting time of visit $v$. Similarly, if a visit $v\in V$ is the last scheduled visit of a caregiver route with qualification $q\in Q_v$ then variable $f_v^q$ is equal to the end time of the visit. If no schedule for a caregiver with qualification $q$ starts (ends) at $v$, variable $k_v^q$ (variable $f_v^q$) is equal to zero. Note that for each qualification type $q \in Q$ there are $|C_q|$ available caregivers, and as a result $k_{v}^{q}$ (resp. $f_{v}^{q}$) can be nonzero for at most $|C_q|$ visits $v \in V$. \cref{tab:notation:MTZ} provides an overview of sets, parameters, and variables used in formulation~\eqref{eq:mtz}.

\renewcommand{\arraystretch}{0.85}
\begin{table}[H]
		\caption{Notation used for MTZ formulation~\eqref{eq:mtz}.}
		\label{tab:notation:MTZ}
		\centering 
		\small
		\begin{tabularx}{\textwidth}{l X} 
	\toprule
	Set & Description  \\
	\midrule
	$V^{\mathrm{o}}$    & Original patient visits \\
	$S$                 & Splittable patient visits ($S \subseteq V^{\mathrm{o}}$) \\
	$V$                 & Potential patient visits ($V^{\mathrm{o}} \cup \{v_p : v \in S, \; p \in \{1,2\}\}$) \\
	$N$                 & Potential patient visits and the artificial start/end of a workday ($V \cup \{0,n+1\}$) \\
	$\Astart$           & Arcs representing the start of a caregivers workday \\  
	$\Acare$            & Arcs representing the performance of a visit and consecutive travel to the next visit \\ 
	$A$                & Allowed travel arcs between patient visits ($\Astart \cup \Acare)$ \\
	$D$                & Visit pairs with temporal dependencies ($D \subseteq V \times V$)  \\
	$Q$                 & Qualification types \\  
	$Q_v$               & Qualification types that satisfy the requirements of visit $v \in N$ \\ 
	$C_q$               & Caregivers with qualification type $q \in Q$ \\
	\midrule
	Parameter & Description \\ 
	\midrule
	$[\alpha_v,\beta_v]$ & Time window during which visit $v \in N$ 
	has to start \\ 
	$d_v$               & Duration of visit $v \in N$ \\ 
	$t_{uv}$             & Travel time from visit $u \in N$ to visit $v \in N$ \\
	$e^q$                & Wage of a caregiver of qualification type $q \in Q$  \\ 
	$\delta_{uv}^{\mathrm{min}},\delta_{uv}^{\mathrm{max}}$ & Min / max starting time difference between $u$ and $v$, $\{u,v\}\in D$, if $u$ starts no later than $v$ \\ 
	\midrule
	Variable & Description \\ 
	\midrule
	$x_a^q\in \{0,1\}$      & Whether a caregiver of type $q \in Q_u \cap Q_v$ uses arc $a=(u,v)$, i.e., successively performs visits $u$ and $v$ \\
	$s_v\in \{0,1\}$           & Whether splittable visit $v \in S$ is split (i.e., if split parts $v_1$ and $v_2$ are performed) \\ 
	$y_v\in \{0,1\}$           & Whether visit $v \in V$ is performed \\ 
	$p_{uv}\in \{0,1\}$        & Whether visit $u$ of visit pair $\{u,v\}\in D$, $u<v$, starts no later than visit $v$\\
	$b_v\ge 0$           & Starting time of visit $v \in V$\\
	$k_v^q\ge 0$         & Starting time of visit $v \in V$ if $v$ is the first visit of a caregiver with qualification $q \in Q_v$\\
	$f_v^q\ge 0$         & End time of visit $v \in V$ if $v$ is the last visit of a caregiver with qualification $q \in Q_v$\\
	\bottomrule
\end{tabularx}
\end{table}

\begin{subequations} \label{eq:mtz}
\begin{small}
    \begin{flalign}
     \mbox{min}\quad & 
 \sum_{v \in V} \sum_{q \in Q_v} e^{q}  (f^{q}_{v} - k^{q}_{v}) \label{eq:mtz:obj} \\
 \text{s.t.}\quad
        & x^q(\delta^{+}(0))\leq  |C_{q}|   &&  \forall q \in Q\label{eq:mtz:caregiver:limit} \\
        & x^q(\delta^{+}(v))  =  x^q(\delta^{-}(v)) && \forall v \in V,\ q \in Q_v \label{eq:mtz:flow:conservation} \\
        & \sum_{q\in Q_v} x^q(\delta^{+}(v))  = y_v && \forall v\in V \label{eq:mtz:inflow}\\
        & y_v =  
       \begin{cases}
           1 & \mbox{if $v\in V^{\mathrm{o}}\setminus S$} \\
           1 - s_v & \mbox{if $v\in S$} \\
           s_u & \mbox{if $\exists u \in S : \ v \in \{u_1,u_2\}$}\\
           \end{cases} && \forall v\in V  \label{eq:mtz:link:y-s}\\        
        &  b_{u} + d_{u} + t_{uv} \sum_{q\in Q_u\cap Q_v} x^q_a \leq   b_{v}  +(\beta_{u} + d_u)(1-\sum_{q\in Q_u\cap Q_v} x^q_a) && \forall a=(u,v)\in A,\ u,v\in V\label{eq:mtz:time} \\    
        & \alpha_v y_v \le b_v \le \beta_v y_v && \forall v\in V \label{eq:mtz:visit:timewindow} \\      
        & b_{v} - \beta_{v}  (1-x_{a}^{q}) \leq   k_{v}^{q} \leq  b_v && \forall v \in V,\ q \in Q_v,\ a=(0,v) \label{eq:mtz:starttimes} \\     %
        &   k_{v}^{q} \leq  (\Tmax-d_{v}) x_{a}^{q}  && \forall v \in V,\ q \in Q_v,\ a=(0,v) \label{eq:mtz:starttimes:forcezero} \\
        & b_{v} + d_v -  (\beta_{v} + d_v)(1-x_{a}^{q})  \leq   f_{v}^{q} \leq  b_v + d_v && \forall v \in V,\ q \in Q_v,\ a=(v,n+1)  \label{eq:mtz:finishtimes} \\
         & f_{v}^{q} \leq  \Tmax x_{a}^{q}  && \forall v \in V,\ q \in Q_v,\ a=(v,n+1)  \label{eq:mtz:finishtimes:forcezero} \\
           &  2 p_{uv} \leq y_u + y_v              && \forall \{u,v\}\in D, u<v \label{eq:mtz:link:p-y}\\              
           & b_v - b_u \le \delta^{\mathrm{max}}_{uv} p_{uv} + \beta_v (2-y_u - y_v)              && \forall \{u,v\}\in D, u<v \label{eq:mtz:synchro:vu:le}\\             
           & b_u - b_v \le \delta^{\mathrm{max}}_{vu}  (y_v-p_{uv}) + \beta_u (2-y_u - y_v)             && \forall \{u,v\}\in D, u<v \label{eq:mtz:synchro:uv:le}\\             
           & b_v - b_u \ge \delta^{\mathrm{min}}_{uv} p_{uv} - \beta_u (2-y_u - y_v) \nonumber  \\ 
        &  \quad \quad \quad \quad - \max\{0,\beta_u-\alpha_v\} (y_u - p_{uv}) &&          \forall \{u,v\}\in D, u<v \label{eq:mtz:synchro:vu:ge}\\              
            & b_u - b_v \ge \delta^{\mathrm{min}}_{vu} (y_v - p_{uv})-(\beta_v+\delta^{\mathrm{min}}_{vu})(2-y_u - y_v)  \nonumber  \\ 
        &  \quad \quad \quad \quad -\max\{0,\beta_v-\alpha_u\} p_{uv} && \forall \{u,v\}\in D, u<v \label{eq:mtz:synchro:uv:ge} \\   
          & y_v\in \{0,1\} && \forall v\in V \\
          & x_a^{q} \in  \{0,1\} && \forall a=(u,v)\in A,\ q\in Q_u\cap Q_v \label{eq:mtz:xvar}\\
         & s_v \in  \{0,1\} && \forall v \in S \label{eq:mtz:svar} \\
         & p_{uv} \in  \{0,1\}  && \forall \{u,v\} \in D, u<v \label{eq:mtz:pvar} \\         
         & b_{v}  \geq  0  && \forall v \in V \label{eq:mtz:bvar} \\
        & k_{v}^{q}\ge  0  && \forall v \in V,\ q \in Q_v \label{eq:mtz:kvar} \\ 
       & f_{v}^{q}   \geq  0  && \forall v \in V,\ q \in Q_v \label{eq:mtz:fvar}
    \end{flalign}
    \end{small}
    \end{subequations}
    
The objective \eqref{eq:mtz:obj} minimizes the operational costs for the HHC provider by multiplying the total working time per qualification type by the wage parameter. 
Constraints~\eqref{eq:mtz:caregiver:limit} ensure that at most the $|C_q|$ available caregivers with qualification type $q\in Q$ can be used. Equations~\eqref{eq:mtz:flow:conservation} are flow conservation constraints for caregiver routes. Equations~\eqref{eq:mtz:inflow} ensure that an appropriately qualified caregiver arrives at the location of a visit if it is performed. Constraints~\eqref{eq:mtz:link:y-s} ensure that each unsplittable visit is performed. For each splittable visit, these constraints ensure that either the original visit or both split parts are performed and that this decision is consistent with the splitting decision. 
Even though visit variables could be easily eliminated using equations~\eqref{eq:mtz:link:y-s}, they are retained as they simplify (the explanation of) the synchronization constraints~\eqref{eq:mtz:link:p-y}--\eqref{eq:mtz:synchro:uv:ge}. Inequalities~\eqref{eq:mtz:time} are commonly used MTZ constraints that ensure that the earliest starting time of a patient visit is equal to the end time of the previous visit plus the travel time between these two visits. These constraints also eliminate subtours. Lower and upper bounds on visit starting times are imposed by \eqref{eq:mtz:visit:timewindow}, which also force the visiting times of visits that are not performed to zero. Together, constraints~\eqref{eq:mtz:caregiver:limit}--\eqref{eq:mtz:visit:timewindow} ensure that all required care is provided by sufficiently qualified caregivers while respecting the time window of each visit.

Constraints~\eqref{eq:mtz:starttimes}--\eqref{eq:mtz:finishtimes:forcezero} ensure the appropriate route starting and finalization times. Therefore, constraints~\eqref{eq:mtz:starttimes} and \eqref{eq:mtz:starttimes:forcezero} ensure that when a visit $v\in V$ is the first visit of a caregiver route of qualification type $q\in Q_v$, variable $k_v^q$ is equal to the starting time $b_v$ of visit $v$ (constraints~\eqref{eq:mtz:starttimes}), while $k_v^q$ equals zero otherwise (constraints~\eqref{eq:mtz:starttimes:forcezero}). Similarly, constraints~\eqref{eq:mtz:finishtimes} and \eqref{eq:mtz:finishtimes:forcezero} ensure that when a visit $v \in V$ is the last visit of a caregiver route of qualification $q \in Q_v$, variable $f_v^q$ is equal to the end time of visit $v$, and equal to zero otherwise.

Constraints~\eqref{eq:mtz:link:p-y}--\eqref{eq:mtz:synchro:uv:ge} ensure that the temporal dependencies for pairs of visits $\{u,v\}\in D$ are met if both visits are performed. Inequalities~\eqref{eq:mtz:link:p-y} ensure that a precedence decision that visit $u\in V$ starts before (or at the same time as) visit $v\in V$ is only possible if both visits are performed. Constraints~\eqref{eq:mtz:synchro:vu:le} and \eqref{eq:mtz:synchro:uv:le} ensure that, depending on the starting order of $u$ and $v$, the second visit is initiated at most $\delta^{\mathrm{max}}_{uv}$ or $\delta^{\mathrm{max}}_{vu}$ time steps after the first one (if both are performed). These constraints are redundant due to their rightmost terms if (at least) one of the two visits is not performed. Similarly, a bound $\delta^{\mathrm{min}}_{uv}$ or  $\delta^{\mathrm{min}}_{vu}$ on the minimum difference between the starting times is imposed by inequalities~\eqref{eq:mtz:synchro:uv:ge} and \eqref{eq:mtz:synchro:vu:ge}, depending on the temporal sequence of tasks $\{u,v\}\in D$ if both are executed. The two right-most terms on the right-hand sides of these constraints make them redundant if one of the visits is not performed and if the left-hand sides of constraints~\eqref{eq:mtz:synchro:vu:ge} or \eqref{eq:mtz:synchro:uv:ge}, respectively, correspond to the opposite temporal starting order.

\subsection{Time-indexed flow formulation}\label{sec:lg}
The time-indexed formulation introduced in this section is defined on a directed time-indexed graph $\LGraph=(\LNodes,\LArcs)$ and assumes that all time-related input parameters (travel times, time windows, etc.) are integer numbers. To this end, we observe that every instance with rational values can be easily transformed into one containing integer numbers by simple rescaling. The node set $\LNodes=\{0\}\cup \LNodesVisit\cup \{n+1\}$ of graph $\LGraph$ consists of start node $0$, end node $n+1$, and the set of nodes $\LNodesVisit=\cup_{v\in V} \LNodesV{v}$ where $\LNodesV{v}=\{v^t\mid t\in \{ \alpha_v, \dots, \beta_v\}\}$ contains one node $v^t$ for each possible initiation time $t$ of visit $v$. Inspired by \citet{FINK2019699}, we classify the set of arcs $\LArcs=\LAstart\cup \LAcare\cup \LAwait$ into starting arcs $\LAstart=\{(0,v^t)\mid v^t\in \LNodesVisit\}$ representing the start of a caregivers workday, waiting arcs $\LAwait=\{(v^t,v^{t+1})\mid v^t,v^{t+1}\in \LNodesVisit\}$ representing one time-unit of waiting at a patient's home, and care arcs $\LAcare=\{(u^t,v^\ell) \mid  u^t,v^\ell\in \LNodesVisit, (u,v)\in A, Q_{u}\cap Q_{v}\ne \emptyset, \ell=\max\{t+d_u+t_{uv},\alpha_v\} \} \cup \{(u^t,n+1) \mid u^t\in \LNodesVisit\}$ $\setminus \{(u^{t},v^{\ell}) \in \LNodesVisit^{2} \mid  \exists w\in S : u,v \in \{w,w_1, w_2\}\}$. Please note that waiting at a patient's home may be required even if a caregiver arrives within that patient's time window due to temporal restrictions with the initiation of other visits. This option is represented by a sequence of waiting arcs. A care arc $(u^t,v^\ell)$ with $u,v\in V$ represents the start of visit $u$ at time $t$ followed by a travel to visit $v\in V\setminus \{u\}$ where the caregiver arrives at time $t+d_u+t_{uv}$. Thus, arc $(u^t,v^\ell)$ where $\ell=\max\{t+d_u+t_{uv},\alpha_v\}$ includes waiting time possibly required at visit $v$ if a caregiver arrives before the earliest time $\alpha_v$ the care at $v$ can start. These care arcs only exist between visit pairs that have a common required qualification type, since a single caregiver must be capable of performing both visits. A care arc $(u^t,n+1)$ indicates that visit $u$ initiated at time $t$ is the last visit of a caregiver's schedule.

Due to the time-indexed nature of graph $\LGraph$, the total working time corresponding to each of its arcs is fixed and we can associate each arc $a\in \LArcs$ with costs 
\begin{equation*}
  c_{a}=\begin{cases}
    0  & \text{if $a = (0,v^{t}) \in \LAstart$,} \\ 
    d_v & \text{if $a = (v^{t}, n+1)\in \LAcare$,} \\     
    \ell - t, & \text{if $a = (u^{t}, v^\ell)\in \LAcare \cup \LAwait$, $v^\ell \ne n+1$,} \\          
  \end{cases} 
\end{equation*}
that correspond to the working time spent when using arc $a$.

Flow formulation~\eqref{eq:ti} uses three sets of binary decision variables. Travel variables $X^q_a\in \{0,1\}$ are equal to one if a caregiver with qualification $q\in Q_u\cap Q_v$ travels along $a=(u^{t},v^{\ell})\in \LArcs$, i.e., performs visit $u$ at starting time $t$ and reaches $v$ at time $\ell$. Split variables $s_v\in \{0,1\}$ indicate whether a splittable visit $v\in S$ is split (in which case split parts $v_1$ and $v_2$ must be performed) or not (in which case visit $v$ must be executed). Finally, for each $\{u,v\}\in D$, $u<v$, precedence variable $p_{uv}\in \{0,1\}$ is equal to one if visit $u$ starts before or at the same time as visit $v$ and equal to zero if visit $v$ starts before visit $u$ or if one of these two visits is not performed at all. \cref{tab:notation:TI} summarizes sets, parameters, and variables used in formulation~\eqref{eq:ti} that are not included in \cref{tab:notation:MTZ}.

\renewcommand{\arraystretch}{0.85}
\begin{table}[H]
		\caption{Additional notation for time-indexed formulation~\eqref{eq:ti}.}
		\label{tab:notation:TI}
		\centering 
		\small
		\begin{tabularx}{\textwidth}{l X}
			\toprule
			Set & Description  \\
			\midrule
			$\LNodesV{v}$       & Time-indexed nodes for visit $v \in V$, each representing a possible initiation time $t$ \\
			$\LNodesVisit$     & Time-indexed nodes for all visits $v \in V$ ($ \cup_{v \in V} \LNodesV{v}$) \\  
			$\LNodes$           & Time-indexed nodes and the artificial start/end of a workday ($\{0\} \cup \LNodesVisit \cup \{n+1\}$). \\
			$\LAstart$          & Arcs representing the start of a caregivers workday at a specific time \\ 
			$\LAcare$           & Arcs representing performing a visit at a specific time and traveling to the next visit \\
			$\LAwait$           & Arcs representing one time-unit of waiting at a patient's home \\
			$\LArcs$           & All arcs between time-indexed nodes ($\LAstart \cup \LAcare \cup \LAwait $). \\
			\midrule
			Parameter & Description\\ 
			\midrule 
			$c_a$                & Cost (working time) associated with arc $a = (u^{t},v^{\ell}) \in \LArcs$ \\ 
			\midrule
			Variable & Description \\ 
			\midrule
			$X_a^q\in \{0,1\}$      & Whether a caregiver of type $q \in Q_u \cap Q_v$ travels along arc $a=(u^{t}, v^{\ell}) \in \LArcs$, i.e., 
			 performs visit $u$ at time $t$ and successively travels to visit $v$ with arrival time $\ell$ \\ 
			\bottomrule
		\end{tabularx}
\end{table}

    \begin{subequations} \label{eq:ti}
    \begin{small}
    \begin{flalign}
    \min\ &  \sum_{a= (u^{t},v^\ell) \in \LArcs} \sum_{q \in Q_{u} \cap Q_{v}} e^{q} c_{a} X_{a}^{q}  \label{eq:ti:obj} \\
             \text{s.t.}\ 
            &  X^q(\delta^{+}(0))\leq |C_{q}|  &&  q \in Q \label{eq:ti:caregiver:limit} \\
            & X^q(\delta^{+}(v^{t})) = X^q(\delta^{-}(v^{t})) && \forall v^{t} \in  \LNodesVisit,\ q \in Q_v \label{eq:ti:flow:conservation} \\
            & \sum_{q\in Q_v} X^{q}(\delta^+(\LNodesV{v})) = 
             \begin{cases}
             1 & \mbox{if $v\in V^{\mathrm{o}} \setminus S$} \\
             1-s_v & \mbox{if $v\in S$}\\
             s_u &\mbox{if $\exists u \in S : \ v \in \{u_1,u_2\}$}\\
            \end{cases} && \forall v\in V \label{eq:ti:visits} \\            
            &  2 p_{uv}  \leq \sum_{q\in Q_u} X^q(\delta^{+}(\LNodesV{u})) + \sum_{q\in Q_v} X^q(\delta^+(\LNodesV{v}))
            && \forall \{u,v\}\in D, \  u<v  && \label{eq:ti:link:p-X} \\
            & p_{uv} + \sum_{q\in Q_u} X^q(\delta^+(u^t)\cap \LAcare) + 
            \sum_{\substack{a=(v^\ell,w^m)\in \Acare : \\ \ell \notin [t+\delta^{\mathrm{min}}_{uv},t+\delta^{\mathrm{max}}_{uv}]}} \sum_{q \in Q_v \cap Q_w} X^q_a \leq 2  &&  \forall u^{t} \in \LNodesV{u}:\{u,v\} \in D,\ u<v \label{eq:ti:syn_12} \\     
            &  \sum_{q\in Q_v} X^q(\delta^{+}(v^t)\cap \LAcare) 
            + \sum_{\substack{a=(u^\ell,w^m)\in \LAcare : \\ \ell\notin [t+\delta^{\mathrm{min}}_{vu},t+\delta^{\mathrm{max}}_{vu}]}} \sum_{q\in Q_u \cap Q_w} X^q_a \leq 1 + p_{uv}
                          && \forall v^{t} \in \LNodesV{v}:\{u,v\} \in D,\ u<v \label{eq:ti:syn_21} \\ 
         & X_{a}^{q} \in \{0,1\} && \forall a=(u,v)\in \LArcs,\ q\in Q_u\cap Q_v \label{eq:ti:Xvar} \\
         & s_v  \in \{0,1\} && \forall v \in S \label{eq:ti:svar} \\
        & p_{uv} \in \{0,1\}  && \forall \{u,v\} \in D, u<v \label{eq:ti:pvar} 
    \end{flalign}
    \end{small}
    \end{subequations}
The objective~\eqref{eq:ti:obj} minimizes operational costs for the HHC provider by multiplying the total working time per qualification type by the wage parameter. Constraints~\eqref{eq:ti:caregiver:limit} ensure for each qualification type that the number of caregiver routes does not exceed the available caregivers. Since $\LGraph$ is acyclic, flow conservation constraints~\eqref{eq:ti:flow:conservation} are sufficient to ensure that each caregiver used completes a route that starts at the artificial depot $0$ and ends at its copy $n+1$. 
Equations~\eqref{eq:ti:visits} ensure that all required care is performed by a sufficiently qualified caregiver. They guarantee that all unsplittable visits are performed, and the execution of either the original visit or both split components for splittable visits. 
Inequalities~\eqref{eq:ti:link:p-X} ensure that a precedence decision indicating that visit $u\in V$ starts before $v\in V$ is only possible if both visits are performed. 
If both visits of a pair $\{u,v\}\in D$ are carried out, constraints~\eqref{eq:ti:syn_12} and \eqref{eq:ti:syn_21} ensure that their temporal dependencies are met. 
To this end, inequalities~\eqref{eq:ti:syn_12} state that if the care of visit $u$ starts at time $t$ and if $u$ is scheduled no later than $v$ (i.e., if $p_{uv}=1$), then $v$ must start during the interval $[t+\delta_{uv}^{\mathrm{min}},t+\delta^{\mathrm{max}}_{uv}]$. This is ensured since none of the arcs indicating a starting time outside of this interval can be used if the first two terms on the left-hand side of a constraint~\eqref{eq:ti:syn_12} are both equal to one. Similarly, constraints~\eqref{eq:ti:syn_21} guarantee that visit $u$ must start in the interval $[t+\delta_{vu}^{\mathrm{min}},t+\delta^{\mathrm{max}}_{vu}]$ if visit $v$ starts at time $t$, both visits are performed, and $u$ starts no earlier than $v$ (i.e., if $p_{uv}=0$).

\section{Solution approach}\label{s: sol_approach}
This section outlines the proposed method to solve the \probabbr. Its design is based on preliminary experiments. Solving the two formulations introduced in the previous section with Gurobi's general-purpose MILP solver yielded the following main insights.
\begin{inparaenum}[(i)]
    \item The time-indexed flow formulation from \cref{sec:lg} performs much better than the MTZ formulation introduced in \cref{sec:mtz}, which suffers from its weak dual bounds.
    \item The large number of variables and constraints in the time-index formulation~\eqref{eq:ti} can present challenges for larger instances.
    \item Finding (good) feasible solutions to the considered instances of the \probabbr is challenging for the general-purpose MILP solver when the size of the instances increases.
\end{inparaenum}

Consequently, we developed a solution method based on the time-index formulation~\eqref{eq:ti} that incorporates pre-processing (see \cref{sec:prepro}) to reduce its size and heuristics to identify high-quality solutions. As detailed in \cref{sec:heuristics}, these heuristics employ either the MTZ or the time-indexed formulation on appropriately selected smaller graphs, and are embedded via callbacks. 

\subsection{Pre-processing}\label{sec:prepro}

We use the following four pre-processing techniques to reduce the size of the time-index graph and, as a consequence, the number of variables in formulation~\eqref{eq:ti}. The first technique only tightens the time windows of the visits and can therefore be applied regardless of the chosen formulation.

\paragraph{Time windows tightening} In accordance with \citet{RASMUSSEN2012698}, we can tighten the time windows of pairs of nodes $\{u,v\}\in D$ that have a predetermined starting order. Assuming, without loss of generality, that $u$ must precede $v$, we can update the time window of $u$ to $[\max \{\alpha_{u}, \alpha_{v} - \delta_{uv}^{max} \}, \min \{\beta_{u}, \beta_{v} - \delta_{uv}^{min} \}]$ and the interval of $v$ to $[\max\{\alpha_{v}, \alpha_{u} + \delta_{uv}^{min}\} ,\min \{\beta_{v}, \beta_{u} + \delta_{uv}^{max}\}]$. Consequently, nodes corresponding to visits $u$ and $v$ outside of these tightened intervals and arcs incident to them can be removed from $\LGraph$. For pairs of nodes $\{u,v\} \in D$ without a predetermined starting order, we can update the time window of $u$ to $[\max\{\alpha_u, \alpha_v -\delta_{uv}^{max}\}, \min \{ \beta_u, \beta_v + \delta_{vu}^{max} \} 
 ]$ and that of $v$ to $[\max\{\alpha_v, \alpha_u -\delta_{vu}^{max}\}, \min \{ \beta_v, \beta_u + \delta_{uv}^{max} \}
 ]$. 
It is important to note that these adjustments can only be made if the execution of one of the visits automatically implies the execution of the other visit (e.g., split components of a visit).

\paragraph{Removal of redundant travel arcs}

Consider two visits $u,v\in V$, $u\ne v$, that can be performed by the same caregiver. Since care arcs possibly include waiting time, graph $\LGraph$ may contain multiple arcs $(u^t,v^{\alpha_v})$ from visit $u$  to the node corresponding to the earliest starting point $\alpha_v$ of visit $v$. Specifically, let $m$ and $\ell$ be the earliest and latest care initiation time of visit $u$, respectively, for which a caregiver can arrive at visit $v$ no later than at time $\alpha_v$. Then, all existing travel arcs between $u$ and $v$ arriving at time $\alpha_v$ are represented by $A'_{uv}=\cup_{t=m}^\ell \{(u^t,v^{\alpha_v})\}$. 

Suppose $u$ does not have synchronization requirements with other visits (i.e., $\nexists w\in V : \{u,w\}\in D$) and a solution exists using a care arc from $A'_{uv}$, such as $(u^t,v^{\alpha_v})$, $m\le t\le \ell$. If $t<\ell$, then for any $t'\in \{t+1, \dots, \ell\}$ a solution with the same objective value is obtained by replacing $(u^t,v^{\alpha_v})$ with $(u^{t'},v^{\alpha_v})$ and waiting arcs $(u^l,u^{l+1})$, $l\in \{t, \dots, t'-1\}$, i.e., the caregiver waits at $u$ instead of at $v$. Hence, all but the last of these travel arcs (i.e., $A_{uv}\setminus \{(u^\ell,v^{\alpha_v})\}$) are removed from $\LGraph$.

For nodes $u$ with synchronization requirements, we also eliminate redundant travel arcs. However, to ensure that we do not eliminate solutions for which no alternative solution (with identical objective value) remains, we adapt the interpretation of the travel variables and the synchronization constraints. The use of arc $(u^t,v^\ell)$ now indicates that the care of $u$ is initialized between the arrival time at node $u$ and time $t$. Consequently, we replace constraints~\eqref{eq:ti:syn_12} by 
\begin{subequations}\label{pre_proc:syn}
\begin{equation}
    p_{uv} + \sum_{q\in Q_u} X^q(\delta^{-}(u^{t}) \setminus \LAwait) + \sum_{\substack{a= (v^{\ell}, w^{m}) \in \LAcare : \\ \ell < t + \delta_{uv}^{\textrm{min}}}} \sum_{q\in Q_{v} \cap Q_{w}} X_a^{q} \leq  2 \label{pre_proc:syn:12a}
\end{equation}
for each $u^{t} \in \LNodesV{u} : \{u,v\} \in D, \ u < v, \ |\delta^{-}(u^{t}) \setminus \LAwait| \geq 1$ and 
\begin{equation}
    p_{uv} + \sum_{q\in Q_u} X^q(\delta^{+}(u^{t}) \cap \LAcare) + \sum_{\substack{a= (w^{m}, v^{\ell}) \in \LArcs \setminus \LAwait : \\ \ell > t + \delta_{uv}^{\textrm{max}}}} \sum_{q\in Q_{w} \cap Q_{v}} X_{a}^{q} \leq  2   \label{pre_proc:syn:12b}
\end{equation}
for all $u^{t} \in \LNodesV{u} : \{u,v\} \in D,\ u < v,\ |\delta^{+}(u^{t}) \setminus \LAwait| \geq 1$. Similarly, we change constraints~\eqref{eq:ti:syn_21} to 
\begin{equation}
      \sum_{q\in Q_v} X^q(\delta^{-}(v^{t}) \setminus \LAwait)  + \sum_{\substack{a= (u^{\ell}, w^{m}) \in \LAcare : \\ \ell < t + \delta_{vu}^{\textrm{min}}}} \sum_{q\in Q_{u} \cap Q_{w}} X_a^{q}
    \leq 1 + p_{uv}     \label{pre_proc:syn:21a}\\
\end{equation}
for every $v^{t} \in \LNodesV{v} :  \{u,v\} \in D,\ u < v,\ |\delta^{-}(v^{t}) \setminus \LAwait| \geq 1$ and 
\begin{equation}
    \sum_{q\in Q_v} X^q(\delta^{+}(v^{t}) \cap \LAcare) + \sum_{\substack{a= (w^{m}, u^{\ell}) \in \LArcs \setminus \LAwait : \\ \ell > t + \delta_{vu}^{\textrm{max}}}} \sum_{q\in Q_{w} \cap Q_{u}} X_{a}^{q}  \leq 1 + p_{uv}     
\label{pre_proc:syn:21b}
\end{equation}
for each $v^{t} \in \LNodesV{v} :  \{u,v\} \in D, \ u < v, \ |\delta^{+}(v^{t}) \setminus \LAwait| \geq 1$.
\end{subequations}

Let $t^{*}$ and $t'$ denote the arrival and departure time at $u$, respectively, and assume that $u$ is performed no later than $v$, i.e., $p_{uv}=1$. Since the care of visit $u$ can be initiated at any time in $[t^{*},t']$, constraints ~\eqref{pre_proc:syn:12a} and \eqref{pre_proc:syn:12b} ensure that $v$ must start in interval $[t^{*}+\delta_{uv}^{\mathrm{min}}, t' + \delta_{uv}^{\mathrm{max}}]$. This is achieved by enforcing an arrival at $v$ before or at $t' + \delta_{uv}^{\mathrm{max}}$, and a departure at or after $t^{*}+\delta_{uv}^{\mathrm{min}}$. 
Similarly, constraints~\eqref{pre_proc:syn:21a} and \eqref{pre_proc:syn:21b} ensure the correct temporal relations if the care at visits $v$ starts before the one at visit $u$.
The replacement of constraint~\eqref{eq:ti:syn_12} and \eqref{eq:ti:syn_21} by \eqref{pre_proc:syn} can also require the creation of additional temporal dependencies. Consider a visit $v \in V$ that has temporal dependencies with $u \in V$ and $w \in V$, i.e. $\{u,v\} \in D$ and $\{v,w\} \in D$. Let $t^{*}$ be the arrival time and $t'$ be the departure time at visit $v$ with $t^{*} \neq t'$. Due to this pre-processing step, the whole interval $[t^{*},t']$ can be used to check the temporal restrictions with $u$ and $w$, whereas previously only time $t'$ was relevant. Hence, it is possible that temporal restriction $\{u,v\}$ utilizes starting time $t^{*}$ at $v$ and temporal restriction $\{v,w\}$ utilizes starting time $t' \neq t^{*}$ at $v$. Since visit $v$ can only start at one moment in time, we specify an additional temporal restriction between $u$ and $w$ to exclude this. For example, when $u$ must start at the same time as $v$, we add an extra temporal restriction between $u$ and $w$ that is similar to the one between $v$ and $w$, i.e., $\delta_{uw}^{\mathrm{min}}= \delta_{vw}^{\mathrm{min}}, \; \delta_{uw}^{\mathrm{max}}= \delta_{vw}^{\mathrm{max}}, \; \delta_{wu}^{\mathrm{min}}= \delta_{wv}^{\mathrm{min}}, \; \delta_{wu}^{\mathrm{max}}= \delta_{wv}^{\mathrm{max}}$. \ref{app:adjust:temp:dep} specifies adjustments for other types of temporal dependencies.

\paragraph{Removal of suboptimal route initialization and ending arcs}
Consider a visit $v\in V$ without synchronization restrictions and assume that its set of time-indexed nodes $\LNodesV{v}$ has at least one node, such as $v^\ell$, for which no outgoing care arc exists, i.e., $\delta^+(v^\ell)\cap \LAcare=\emptyset$. In this case, an optimal solution of the \probabbr cannot include a caregiver route starting with visit $v$ at time $\ell$ unless $v$ is the only visit scheduled for this caregiver. This holds true since it is cheaper to start the route later and avoid the waiting time required after arriving at visit $v$ at time $l$ (since $v^\ell$ has no outgoing care arc). Similarly, there is always an alternative optimal solution in the problem variant \probabbrtt that avoids this required waiting time. Thus, we remove all starting arcs $(0,v^\ell)$ to nodes $v^\ell$ without synchronization requirements without outgoing care arcs. Similar arguments can be used to observe that a visit $v$ without synchronization requirements initiated at time $t$ cannot be the last visit of an optimal solution of the \probabbr (an alternative optimal solution will also exist in variant \probabbrtt) if the corresponding node has no ingoing care arc and if the caregiver performs at least two visits. Thus, we remove all arcs $(v^\ell,n+1)$ for these nodes. 
Observe that these preprocessing steps may eliminate the option of performing a route with a single visit only. Finally, we check if at least one node exists from $\LNodesV{v}$ incident to both a starting and an ending arc. If no such node exits, we pick one to which we add the required arcs. To select this node, we search in the following order: i) the earliest node with an outgoing arc, ii) the earliest node with an ingoing arc, iii) the earliest node of the visit. 

\paragraph{Replacement of consecutive waiting arcs}
The previous pre-processing step may lead to sequences of degree-two nodes $(v^t, \dots v^\ell)$, $\ell\ge t$, corresponding to the same visit $v\in V$ that are incident to one ingoing and one outgoing waiting arc, i.e., $\delta^-(v^m)=\{(v^{m-1},v^m)\}$ and $\delta^+(v^m)=\{(v^{m},v^{m+1})\}$ for all $m\in [t, \dots, \ell]$. 
If such a sequence cannot be extended with another node, it has maximal length and we replace nodes $\{v^t, \dots, v^\ell\}$ and their incident waiting arcs by a single waiting arc $(v^{t-1},v^{\ell+1})$.

\subsection{Adding heuristic components}\label{sec:heuristics}
After applying the pre-processing techniques, we solve the time-indexed formulation \eqref{eq:ti} using a general-purpose MILP solver enriched with tailored primal and improvement heuristics. The development of these heuristic procedures was triggered by the fact that preliminary experiments showed that the solver has difficulty identifying (high-quality) solutions, in particular for larger-sized instances, and can spend a long time improving the initial solutions it found. 
We therefore developed two versions of a primal heuristic that aim to identify a first feasible solution. They are applied at every node of the branching tree if the associated LP solution has fractional values and if no feasible solution has been found yet. The main difference between the two versions is that the first one is based on the time-indexed formulation~\eqref{eq:ti}, while the second one builds upon the MTZ formulation~\eqref{eq:mtz}. In addition to the primal heuristic, we apply an improvement heuristic that focuses on improving the timing of caregiver routes of integer solutions. Both the primal and improvement heuristics are ILP-based and implemented using callbacks of Gurobi's general-purpose solver, but any preferred general-purpose solver can be used. To limit the total time spent on these callbacks, we apply a time limit corresponding to a certain proportion of the total amount of available time. 
A detailed description of these heuristics is given in the following paragraphs, where we assume that a bar $\bar{\cdot}$ over (a set of) variables indicates their values in a given (possibly fractional) candidate solution. A graphical overview of the solution approach is presented in \cref{fig:flow_chart_sol_appr}.

\begin{figure}[h]
\centering
\includegraphics[width=1.0\linewidth]{graphs/MethodeFlowChart_upd_3ptline.png}
\caption{Flowchart of the solution approach.}
\label{fig:flow_chart_sol_appr}
\end{figure}

\paragraph{Time-indexed based primal heuristic}
This heuristic is based on solving the time-indexed formulation~\eqref{eq:ti} on a graph $\bar{\LGraph}=(\bar{\LNodes},\bar{\LArcs})$, $\bar{\LNodes} \subseteq \LNodes$ that only uses a subset of the nodes of graph $\LGraph$. Graph $\bar{\LGraph}$ is initialized as the subgraph induced by all arcs $\bar{\LArcs}$ used in the solution of the current LP relaxation, i.e., $\bar{\LArcs}=\{a=(u^t, v^\ell)\in \LArcs : \exists q\in Q_u\cap Q_v : \bar{X}_a^q>0\}$, where $\bar{\LNodes}$ denotes the set of nodes induced by $\bar{\LArcs}$. We further add a set of arcs to increase the probability of finding a feasible solution. In particular, for each visit $v\in V$ for which at least one node $v^t$ is 
in $\bar{\LNodes}$ we add the following arcs:
\begin{itemize}
    \item a route start arc $(0,v^\ell)$ to the earliest copy of visit $v$ included in $\bar{\LNodes}$ for which an outgoing care arc is included in $\bar{\LArcs}$;
    \item a route ending arc $(v^\ell,0)$ from the latest copy of visit $v$ included in $\bar{\LNodes}$ for which an ingoing care arc is included in $\bar{\LArcs}$; and
    \item waiting arcs $(v^\ell,v^m)$ for each pair of copies of visit $v$ within $\bar{\LNodes}$ that have consecutive times.
\end{itemize}

For each visit $v\in V$, we possibly add arcs ensuring the existence of at least one option for performing a route with a single visit only. If no such node exists, then we search for a node on which we ensure this option in the following order: the earliest node with an outgoing arc, followed by the earliest node with an ingoing arc, and finally, the earliest node of the visit. Lastly, for visits $v\in V$ with synchronization requirements we also add a route start arc $(0,v^\ell)$ and a route ending arc $(v^m,n+1)$ to the earliest and latest copy of visit $v$ included in $\bar{\LNodes}$, respectively. 

\paragraph{MTZ-based primal heuristic}
This variant of our primal heuristic consists of two phases. In the first phase, we apply the MTZ formulation~\eqref{eq:mtz} to a subgraph $\bar{G}$ of $G$ induced by the current LP relaxation, i.e., the graph induced by arc set $\bar{A}=\{(u,v) : \exists a=(u^t, v^\ell) \in \LArcs, q\in Q_u\cap Q_v : \bar{X}_a^q>0\}$. The aim of this phase is to either quickly identify a feasible solution or show that graph $\bar{G}$ does not contain any solution. Consequently, we stop the solution process when a first solution is found. 
If a solution is found in the first phase, then we apply a second phase that aims to improve this solution. Here, for each visit $v\in V$ we add a starting arc $(0,v)$ and ending arc $(v,n+1)$ to graph $\bar{G}$ in case they are not used in the solution of the current LP relaxation. Afterwards, we solve MTZ formulation~\eqref{eq:mtz} for this extended graph.
 
\paragraph{Improvement heuristic}
The timing of caregiver routes may be suboptimal for an incumbent solution identified by the solver or the solution identified using the time-indexed approach described above. Thus, for each such solution we apply the following ILP-based improvement heuristic to optimize the visit starting times. We construct a graph $\bar{G}=(N,\bar{A})$ whose set of arcs consists of route start and end arcs $\cup_{v\in V} \{(0,v),(v,n+1)\}$ and all arcs $\{(u,v) : \exists a=(u^t, v^\ell) \in \LArcs, q\in Q_u\cap Q_v : \bar{X}_a^q=1\}$ corresponding to travel options used in the solution triggering the application of the improvement heuristic. Given the small number of travel options, we simply solve the MTZ formulation~\eqref{eq:mtz} on subgraph $\bar{G}$ and return the solution as a potential new incumbent to the solver.

\section{Benchmark instances}\label{sec:instances}
To explore the consequences of task-splitting for home healthcare planning, we generated a set of benchmark instances based on those introduced by \citet{Bredstrom2008}. These instances were designed to mimic an HHCRSP with strict synchronization and have also been utilized in computational experiments of other papers that study HHCRSPs with synchronization (see, e.g., \cite{Masmoudi2023, Euchi2020, Decerle2019,Haddadene2019}).
We adapt the instances of \cite{Bredstrom2008} with a small visiting time window (five time window widths are available) which include one caregiver for every five visits, a planning horizon of 9 hours\footnote{In 4 instances, the planning horizon could potentially extend until the original latest starting time of an unsplit visit, resulting in an increase of the planning horizon by 2.0\%, 2.7\%, 0.1\%, and 7.0\%.}, and strict synchronization between one pair of visits out of every ten patient visits. We introduce qualification requirements for visits, qualifications of caregivers, and task-splitting possibilities. Three hierarchically ordered qualification levels are considered (1,2, and 3, where level 3 caregivers can perform all care). Similar to \cite{Yin2023} and \cite{Qiu2021} we assign costs $e^1= 1$, $e^2= 2$, and $e^3= 3$ to them. Generating the qualifications ourselves instead of using, e.g., the randomly generated services used in \cite{AitHaddene2016} provides us the flexibility to consistently introduce task-splitting possibilities and combine them with different scenarios (i.e., varying staff availability and visit requirements) in each of the instances.

For each original \cite{Bredstrom2008} instance, we create different scenarios that vary the qualification requirements of the visits and qualifications of the available staff. Specifically, we vary the proportion of visits that require specific qualification levels using the following three profiles: 
\begin{inparaenum}[(i)]
    \item general support requirements (50\% level 1, 25\% level 2, 25\% level 3) 
    \item balanced requirements (33.3\% level 1, 33.3\% level 2, 33.3\% level 3) 
    \item medical requirements (25\% level 1, 25\% level 2, 50\% level 3). 
\end{inparaenum} 
We generate the required qualification type of each visit independently from its duration and maintain the required qualification level for as many visits as possible when transitioning between scenarios. We also vary the qualification levels of the caregivers, considering four scenarios: 
\begin{inparaenum}[(i)]
    \item practically trained (50\% level 1, 25\% level 2, 25\% level 3);
    \item moderately trained (25\% level 1, 50\% level 2, 25\% level 3);
    \item medically trained (25\% level 1, 25\% level 2, 50\% level 3); and
    \item only medically trained (0\% level 1, 0\% level 2, 100\% level 3). 
\end{inparaenum}
We derive ten benchmark instances for each initial instance by \citet{Bredstrom2008}. These instances are the result of combining the three different visit qualification profiles with the first three staff compositions plus the one in which qualification levels become irrelevant (100\% level 3 staff).

\citet{Bredstrom2008} consider 5 instances of 20 visits, 3 instances of 50 visits, and 2 instances of 80 visits. We scale down the 50 and 80 visit instances of \cite{Bredstrom2008} by selecting a random subset of visits while maintaining the ratio of synchronization visits. From each of those instances, we create new instances for the \probabbr of 20, 30 and 40 visits. For these generated instances, the visits of the smaller instances are included in the larger ones, and each visit has the same characteristics and required qualification level within a certain scenario.

The variation in visit requirements and staff composition together with the variation in instance size results in a variety of instances to challenge the proposed solution approaches. It also allows us to gain insights into the potential impact of task-splitting on a highly diverse test set. 

To complete the benchmark instances, we introduce the possibility of splitting tasks for visits of at least one hour. \cref{tab:splittablevisits:overview} provides insight into the number of splittable visits within the instances considered. The scheduling characteristics of the split components correspond to those of the unsplit visits (combined duration, qualification requirements, permitted starting and completion time). However, for visits with a duration of 75 minuter or more, the combined \emph{duration} of both split components can be reduced or increased by 15 minutes with probability $\frac{1}{3}$. The duration of each split component is determined by randomly splitting the duration of the original visit while ensuring a minimum duration of 30 minutes for each of them. We also relax the required \emph{qualification level} of one randomly selected split component to level 1 with probability 0.75. In addition, we increase the latest starting time of one of the split components
by one hour with probability 0.75. Finally, we randomly impose one of the following three \emph{temporal dependencies} between the split components uniformly at random: no temporal dependency, precedence, and disjunction. 
For disjunctive and precedence temporal dependencies, we ensure that waiting time potentially removed from a caregiver visit due to a split is still present for the patient. 
Furthermore, two splittable visits that require strict synchronization are split in the same way (except possible caregiver qualification relaxation), and their split parts must be performed together (precedence or disjunction is enforced between two split parts). 

\renewcommand{\arraystretch}{0.65}
\begin{table}
\centering
	\caption{Key characteristics of the generated instances.}\label{tab:splittablevisits:overview}
\begin{tabular}{c c c c c c}
	\toprule
    \multirow{2}{*}{\shortstack{number of \\ original visits}} & \multirow{2}{*}{\shortstack{number of \\ caregivers}} & \multirow{2}{*}{\shortstack{number of \\ instances}} & \multicolumn{3}{c}{number of splittable visits} \\ 
     \cmidrule(lr){4-6}	
    & & & average & minimum & maximum \\ 
    \midrule  
	20 & 4 & 10 & 10.8 & 7 & 16  \\
        30 & 6 & 5 & 17.2 & 12 & 22   \\ 
        40 & 8 & 5 & 22.8 & 19 & 28 \\
	\bottomrule
\end{tabular}
\end{table}

\section{Computational results}\label{sec:results}
In this section, we report and discuss the results of our computational experiments performed on the instances described in the previous section. 
This section is structured according to the following two main goals of our computational study:
\begin{enumerate}
    \item \emph{Performance analysis:} \cref{sec:results:performance} studies the performance of the suggested solution methods and the extent to which the novel aspect of task-splitting contributes to the difficulty of solving instances of the \probabbr. 
    \item \emph{Impact of integrating task-splitting:} \cref{sec:results:impact} provides insights into the potential impact of optimizing task-splitting decisions for HHC providers and analyzes which conditions and scenarios influence this impact.
\end{enumerate}
Our solution method was implemented in the programming language \textit{Julia} using Gurobi 10.0.0. 
All experiments were conducted on a single core of an \textit{AMD Genoa 9654} processor of the Dutch National Supercomputer Snellius. The time allotted for solving a single instance was limited to a maximum of five hours, while the cumulative time used by heuristics that are incorporated via callbacks (see, \cref{sec:heuristics}) was limited to a maximum of 30 minutes. 

We note that the combination of high qualification requirements for visits (i.e., medical requirements) with a practically or moderately trained staff renders some instances infeasible. 
In the following section, we will therefore consider only the subset of 159 feasible instances out of a total of 200 created (75/100 with 20 visits, 44/50 with 30 visits, and 40/50 with 40 visits). 

In the computational results, the following abbreviations will be used to denote different configurations of our solution method: \TI refers to solving the time-indexed formulation using Gurobi. \TIHTI and \TIHMTZ refer to the two variants, which also include the heuristics introduced in \cref{sec:heuristics}, where the primal heuristic is used in its time-indexed or MTZ-based variant, respectively. 
A superscript minus (e.g., \TINOSPLIT) will be used for the variant of \probabbr without task-splitting (i.e., $s_v = 0$, $\forall v \in S$) and a superscript plus (e.g., \TIALLSPLIT) when solving the variant of \probabbr in which all splittable visits are split (i.e., $s_v = 1$, $\forall v \in S$). The latter \probabbr variant (where all splittable visits must be executed as split visits) represents an extreme approach of full task-splitting, where all splitting decisions are enforced in advance and are not optimized based on their impact on the overall caregiver planning. 
The pre-processing techniques described in \cref{sec:prepro} are used in all variants and in the order in which they are presented, since this sequence results in the highest reduction of the number of arcs in the time-indexed graph. Throughout the following sections, the number of instances of a certain type is referred to by $\#$ and the term \emph{size} is used to indicate the number of original visits (i.e., 20, 30, or 40) within an individual instance. 
Furthermore, the setting in which all splittable visits need to be split will be referred to as mandatory splitting. 
The term \emph{task-splitting} refers to the general case when splitting visits is optional.

\subsection{Performance analysis}\label{sec:results:performance}

The introduction of task-splitting possibilities is expected to have a large impact on the time required to solve an instance. Therefore, this aspect is first discussed for the \TI method. Afterwards, insight is provided into the performance of the two types of proposed solution approaches, in order to evaluate the effectiveness of enriching the solver with either of the two primal heuristics in combination with the improvement heuristic.

\paragraph{Difficulty of the \probabbr}

\begin{figure}
\centering
\includegraphics[width=1.0\linewidth]{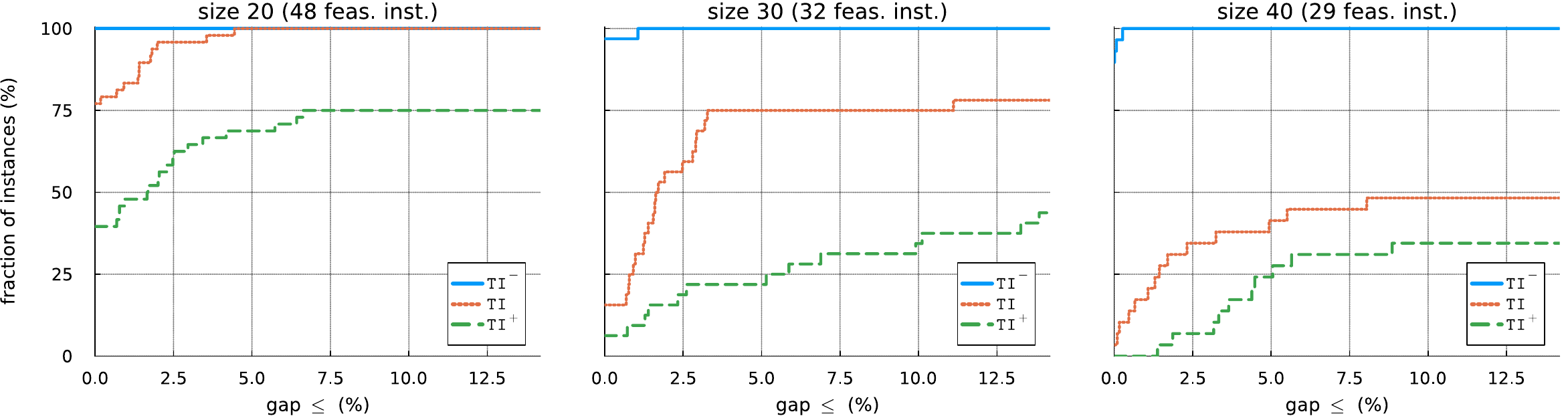}
\caption{Optimality gaps of \TINOSPLIT, \TI, and \TIALLSPLIT on the subset of instances feasible without task-splitting.}
\label{fig:perf_plain_solver}
\end{figure}

The possibility to split longer patient visits into two separate visits increases the size of an instance. For the considered instances, on average, 56\% of the visits are splittable. Consequently, the average number of patient visits included in an instance is more than doubled, while the number of (travel) decision variables after pre-processing and presolving is approximately six times as high. This large increase in the number of decision variables is due to the relaxed scheduling conditions of the split part.
In the following section, we analyze the impact of task-splitting on the difficulty of solving \probabbr instances with Gurobi's general-purpose solver (i.e., method \TI).
\cref{fig:perf_plain_solver} compares the relative amounts of instances with an optimality gap below a certain threshold after five hours without task-splitting, with task-splitting decisions, and when all splittable visits need to be split, differentiating between the size of the instances. 
If the plot for a particular method does not reach 100\%, then this signals that no feasible solution was found for the remaining fraction of the instances. The x-axis includes the maximum optimality gap among all instances for that method. 
In this paragraph, we restrict our attention to the subset of 109 instances that are feasible without task-splitting. 
The results clearly illustrate the additional computational challenge posed by task-splitting. Without task-splitting, 105 of the 109 feasible instances are solved optimally within the given time limit. Moreover, the largest optimality gap for the four remaining instances is only 1.06\%. 
With task-splitting possibilities, only 43 instances could be solved optimally, no feasible solution could be found for 22 instances, and the optimality gaps are significantly larger (up to 11.2\%) for the remaining 44 instances. The results also indicate that the variant with mandatory splitting (i.e., solving \TIALLSPLIT) seems to be the most difficult one from a computational perspective. This variant enforces task-splitting in advance and therefore excludes the optimization of task-splitting decisions based on the eventual effect on caregiver routes and schedules. In this case, only 21 instances are solved optimally and feasible solutions (with gaps up to 13.8\%) are found for 39 instances. Additional experiments with longer time limits (cf.\ \cref{sec:results:impact}) revealed that only 103 of the 109 instances feasible without task-splitting are feasible when all splittable visits need to be split. Thus, no feasible solution could be identified within 5 hours for 43 feasible instances when using $\TIALLSPLIT$. 
These results clearly demonstrate that the \probabbr is much more challenging to solve than its simpler counterpart without task-splitting and that splitting all splittable visits in advance, instead of optimizing task-splitting decisions, seems to be the most difficult option. The additional complexity of task-splitting is specifically observed for the larger instances. The fact that the solver does not find a feasible solution for a large proportion of the instances clearly motivates the development of heuristics and their integration into the exact solution method.

\paragraph{Impact of heuristics} 

\begin{figure}
\centering
\includegraphics[width=1.0\linewidth]{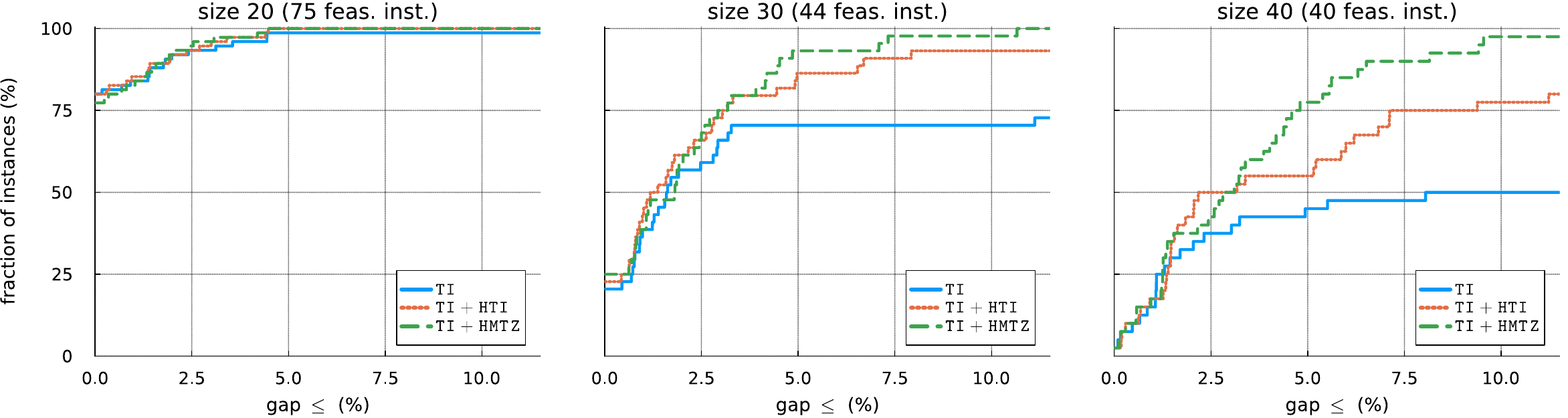}
\caption{Optimality gaps of \TI, \TIHTI, and \TIHMTZ.}	
\label{fig:gaps_3_approaches}
\end{figure}

Motivated by the preceding results, we proceed with the analysis of the impact of the heuristics introduced in \cref{sec:heuristics}. To this end, \cref{fig:gaps_3_approaches} shows the cumulative percentages of instances for which the optimality gap is below a certain threshold when the time-indexed formulation (\TI) is solved in isolation or when it is solved in conjunction with the time-indexed or the MTZ-based heuristic (\TIHTI and \TIHMTZ). In these experiments, all 159 instances that are feasible for the \probabbr are considered. If the plot for a particular method does not reach 100\%, then this signals that no feasible solution was found for the remaining fraction of the instances. The x-axis includes the maximum optimality gap among all instances for that method. 
\cref{fig:gaps_3_approaches} shows that both variants of the heuristic result in significant performance improvements, especially for instances of size 30 and 40.
For small instances of size 20, we observe that they can be solved relatively well without the use of heuristics. Consequently, the observed differences between the three variants are minimal for these instances. 
Both versions with heuristics identify a feasible solution for one additional instance, and the optimality gaps are nearly identical.
However, a clear performance gain is observed for the instances of size 30, for which \TI could not find any feasible solution in 12 cases. Here, \TIHTI finds feasible solutions for 9 additional instances, while feasible solutions to all instances are found by the MTZ variant \TIHMTZ. 
The performance gain becomes even larger for the instances of size 40. In these instances, \TI, \TIHTI, and \TIHMTZ find feasible solutions for 20, 32, and 39 out of the 40 instances, respectively. 
In addition to finding more feasible solutions, we also observe that including the heuristics significantly reduces the final optimality gaps in instances of size 30 and size 40. We observe that variant \TIHMTZ with the MTZ-based primal heuristic outperforms the time-indexed variant \TIHTI with respect to the number of instances for which solutions are found and when considering the final optimality gaps.
The positive impact of the heuristics, the clear benefits of \TIHMTZ over \TI and \TIHTI, as well as the significant contributions of the improvement step are also confirmed by the more detailed results discussed in \ref{app:perf_heuristics}. Overall, we conclude that variant \TIHMTZ outperforms the variants \TIHTI and \TI, and the obtained results demonstrate the performance of the proposed heuristic procedure.

\subsection{Impact of integrating task-splitting}\label{sec:results:impact}

The impact of introducing and optimizing task-splitting decisions for HHC providers is discussed in four parts. First, we explore the effect on staff requirements. Second, we discuss the consequences for the work schedules of caregivers. Third, we describe the change in operational costs achieved through the utilization of task-splitting options. Fourth, we evaluate the percentage of split options utilized in different scenarios. Since (proven) optimal solutions are not available for most task-splitting instances, the analyses compare the best available solution without task-splitting to the best available solution when task-splitting is permitted. We emphasize that although we carefully selected reasonable parameter values in our benchmark instances they are still artificial ones. Therefore, the goal of this analysis is to explore the potential benefits of task-splitting in different scenarios.

\paragraph{Required staff and staff composition}
One primary motivation for considering task-splitting is that it may enable a better use of the available staff. To ascertain under which conditions such an effect can be observed, we present in \cref{tab:number_of_feasible_instances} the number of generated instances for which a feasible schedule for caregivers exists without task-splitting, with task-splitting decisions, and when all splittable visits need to be split. 
We first observe that the additional planning flexibility resulting from task-splitting significantly increases the number of cases for which all patient visits can be performed with the available caregivers. For instances of size 20, task-splitting increases the number of feasible instances from 48 to 75, representing a +56\% increase. For instances of size 30, the increase is from 32 to 44 (+38\%) and for instances of size 40, the increase is from 29 to 40 (+38\%). We also observe that forcing all splittable visits to be split, instead of only allowing it, decreases the number of feasible instances of size 20 and 30 by 20 (-27\%) and 2 (-5\%), respectively. 
This indicates the potential benefits of integrating task-splitting, scheduling, and routing decisions for HHC providers and patients over a sequential approach in which task-splitting decisions are made upfront.

The results for the benchmark instances further indicate that task-splitting allows for the performance of patient visits with a less medically trained staff compared to a situation without task-splitting possibilities. This can be observed from the results for scenarios with balanced patient requirements in combination with a practically or moderately trained staff, where most infeasible instances become feasible thanks to task-splitting. This is also noticed in case of mandatory task-splitting (for the instances of size 30 and 40 and partly for instances of size 20) where no optimization of task-splitting decisions is performed. 
In contrast, most instances remain infeasible for a practically or moderately trained staff in scenarios with medical patient visit requirements. 
In addition to the extra planning options offered by the ability to split parts between different caregivers at different moments in time, other advantages of task-splitting are also utilized. 
The amount of tasks within a certain time frame that require more medically trained caregivers can be reduced by performing one split part of a splittable visit by a more practically trained caregiver, if allowed. In instances that become feasible through task-splitting, this advantage is utilized in approximately 73\% of the split visits. 
Task-splitting can reduce the number of tasks that must be performed within a certain timeframe since some split visits have an extended, relaxed starting time window (see \cref{sec:instances}). However, for the test instances, this factor seems to play a less prominent role than the qualification relaxation and is exploited in $\pm$35\% of the splits in instances that became feasible through task-splitting.

We conclude that task-splitting can help increase planning options and reduce the number of tasks that require a specifically trained caregiver type. Thereby, it is advantageous to integrate task-splitting, routing and scheduling decisions instead of deciding how to organize care tasks in advance. The larger planning flexibility ensures that tasks are only split if this is beneficial for the overall outcome. Integrating task-splitting decisions can potentially help to counteract issues arising from staff shortages, especially for medically trained caregivers.

\begin{table}
\centering
	\caption{Number of feasible instances without task-splitting (no split), with task-splitting possibilities (split), and with mandatory splitting of all splittable visits (full split).} 
    \label{tab:number_of_feasible_instances}
    \small
\resizebox{\columnwidth}{!}{
\setlength\tabcolsep{1mm}
\begin{tabular}{rrrcccccccccccc}
	\toprule
& & & \multicolumn{12}{c}{number of feasible instances} \\ 
\cmidrule(lr){4-15}
& & & \multicolumn{3}{c}{prac. train. staff} & \multicolumn{3}{c}{moder. train. staff} & \multicolumn{3}{c}{medic. train. staff} & \multicolumn{3}{c}{only medic. staff} \\
\cmidrule(lr){4-6} \cmidrule(lr){7-9} \cmidrule(lr){10-12} \cmidrule(lr){13-15}
size & \# & visit req. & no split  &  split & full split & no split  &  split & full split & no split  &  split & full split & no split  &  split & full split  \\ 
\midrule
\multirow{3}{*}{20} & \multirow{3}{*}{10}& general & 4 & 10 & 6 & 5 & 10 & 7 & 10 & 10 & 8 &        &        &  \\
& & balanced & 1 & 5  & 4 & 3 & 10 & 6 & 8  & 10 & 8 & 10     & 10     & 8  \\
& &  medical& 0 & 0  & 0 & 0 & 0  & 0 & 7  & 10 & 8 &        &        &    \\ \midrule 
\multirow{3}{*}{30} & \multirow{3}{*}{5} & general & 5 & 5  & 5 & 5 & 5  & 5 & 5  & 5  & 5 &        &        &    \\
& & balanced& 0 & 5  & 5 & 3 & 5  & 5 & 5  & 5  & 5 & 5      & 5      & 5\\  
& & medical &0 & 1  & 1 & 0 & 3  & 1 & 4  & 5  & 5 &        &        &   \\ \midrule  
\multirow{3}{*}{40}&  \multirow{3}{*}{5} & general & 5 & 5  & 5 & 5 & 5  & 5 & 5  & 5  & 5 &        &          \\
& & balanced& 0 & 5  & 5 & 0 & 5  & 5 & 5  & 5  & 5 & 5      & 5      & 5   \\
& & medical& 0 & 0  & 0 & 0 & 0  & 0 & 4  & 5  & 5 &        &        &   \\
	\bottomrule
\end{tabular}
}
\end{table}

\paragraph{Staff schedules}
To explore the impact of task-splitting on the daily work schedule of caregivers, \cref{tab:percentage_patient_time_minimization_operat_costs} lists the average percentages of the total working time that caregivers spend on performing care. Only instances that are feasible without task-splitting are considered in this analysis. We observe that concerns may arise regarding the impact of task-splitting on caregiver work satisfaction. They include increased travel and/or waiting time or overly dense plannings. However, for the instances considered, we observe that on average, prior to splitting, 84.7\% of the total working time of the caregivers consisted of patient visits, while after splitting, this decreased slightly to 84.0\% (these numbers can be obtained as the weighted averages of the results given in \cref{tab:percentage_patient_time_minimization_operat_costs}). Therefore, we conclude that optimizing task-splitting does not seem to significantly change the percentage of working time caregivers spend providing care for the considered instances. This could be of interest to HHC providers when counteracting potential staff concerns. 
Furthermore, \cref{tab:percentage_patient_time_minimization_operat_costs} does not demonstrate any clear relationships between the change in patient time and the considered scenario (e.g., see general support visit requirements in combination with different staff compositions or a medically trained staff composition with different visit requirements for size 30 instances). The results for the benchmark instances indicate that a higher percentage of time is spent on patient care in scenarios where more medically trained caregivers become available (compare combination of general support visit requirements with a practically or only medically trained staff) or where visits require less medical tasks (moving from medical visit requirements to general support requirement for a medically trained staff). 

Our findings on the test instances indicate that task-splitting has no significant impact on the total working time of the caregivers. The total working time is, on average, reduced by only 0.8\% across the instances considered. Indeed, the relative increase in travel time due to task-splitting (+9\%) is compensated by a reduction in patient visiting and waiting time. 

Although the impact of task-splitting on the total working time is limited for the instances considered, it does impact the division of work among caregivers. In the set of test instances, the average percentage of total working time assigned to level 3 caregivers decreased by 4.3 percentage points, indicating a reduction in their participation in tasks requiring a level 1 or level 2 caregiver. 
The consequence of level 3 caregivers performing fewer care tasks is that mainly level 1 caregivers work longer.
Therefore, one of the potential benefits of task-splitting may be that caregivers spend relatively more time on tasks aligned with their education. An overview of the division of working time among the three qualification levels with and without splitting is presented in \ref{app:devision_working_time}.

\begin{table}
\centering
\caption{Average percentage of the total working time spent on patient care.}
\label{tab:percentage_patient_time_minimization_operat_costs}
\small
\begin{tabular}{rrrcccccccc}
	\toprule
& & & \multicolumn{8}{c}{average percentage of working time spent on patient care (\%)} \\
\cmidrule(lr){4-11}
& & & \multicolumn{2}{c}{prac. train. staff} & \multicolumn{2}{c}{moder. train. staff} & \multicolumn{2}{c}{medic. train. staff} & \multicolumn{2}{c}{only medic. staff} \\
\cmidrule(lr){4-5} \cmidrule(lr){6-7} \cmidrule(lr){8-9} \cmidrule(lr){10-11}
size & $\#$ & visit req. & no split & split & no split & split & no split & split & no split & split \\ 
	\midrule
\multirow{3}{*}{20} & \multirow{3}{*}{10} & general         & 78.0 & 77.2 & 80.8 & 80.5 & 82.6 & 82.5 &      & \\
                                                & & balanced& 80.3 & 82.0 & 77.2 & 79.0 & 80.4 & 80.8 & 86.2 & 85.8   \\
                                                & &  medical&      &      &      &      & 79.3 & 78.5 &      & \\   \midrule 
\multirow{3}{*}{30} & \multirow{3}{*}{5}  & general         & 85.2 & 83.9 & 87.0 & 85.8 & 87.5 & 86.5 &      &  \\ 
&                                                 & balanced&      &      & 84.2 & 83.1 & 86.1 & 85.2 & 89.8 & 89.1    \\
                                                & & medical &      &      &      &      & 85.1 & 83.2 &      &      \\ \midrule  
\multirow{3}{*}{40}&  \multirow{3}{*}{5} & general          &86.2 & 84.0 & 87.8 & 86.4 & 88.5 & 87.9 &      &\\
                                                & & balanced&     &      &      &      & 87.3 & 86.1 & 91.3 & 90.4    \\
                                                 & & medical&      &      &      &      & 85.6 & 84.0 &      & \\
	\bottomrule
\end{tabular}
\end{table}
\paragraph{Operational costs} 
The possibility of task-splitting not only reduces the staffing requirements for the instances studied, but also allows for schedules with lower operational costs. This is achieved by performing only those splits that contribute to a better overall planning. 
\cref{tab:decrease_in_oper_cost} summarizes the average decrease in operational costs due to task-splitting in columns four to seven. 
The smallest gains occur in the instances where each caregiver can perform all patient visits. In the scenario with only medically trained staff, the flexibility gained through task-splitting reduces the total costs by more than 1.1\%. 
This decrease is primarily achieved through splits that reduce the combined patient time. The gains from task-splitting increase for scenarios that differentiate between the qualifications of the caregivers. This effect is particularly strong in the generated scenarios where a large proportion of the caregivers are unable to perform a subset of the visits. 
Focusing, for instance, on the general support scenario, we observe that the average gain due to task-splitting increases each time the staff composition shifts more to a practically trained staff. 
Without task-splitting, medically trained caregivers must perform all visits that include a highly medical task. However, due to splitting, all types of caregivers (i.e., 1, 2 and 3) can perform split parts that only require a practically trained caregiver.

\begin{table}
	\centering
	\caption{Average operational cost decrease due to task-splitting and average percentage of utilized splits for practically trained (prac), moderately trained (moder), medically trained (med), and only medical staff (only med) staff profiles. 
	}
	\label{tab:decrease_in_oper_cost}
	\small
		\begin{tabular}{rrrcccccccc}
			\toprule 
			& & & \multicolumn{4}{c}{average decrease in operational cost (\%)} & \multicolumn{4}{c}{average percentage of utilized splits (\%)} \\
		\cmidrule(lr){4-7} \cmidrule(lr){8-11}
		size & $\#$ & visit req. & prac & moder & med & only med & prac & moder & med & only med \\
		\midrule
		\multirow{3}{*}{20} & \multirow{3}{*}{10} & general         &3.83 & 2.39 & 4.33 &      &40 & 29 & 29 &    \\
		& & balanced&4.68 & 4.26 & 5.97 & 1.16 & 43 & 35 & 32 & 21 \\
		& &  medical&     &      & 8.08 &   &   &    & 47 &  \\ \midrule 
		\multirow{3}{*}{30} & \multirow{3}{*}{5}  & general         &7.26 & 3.37 & 2.23 &  &     44 & 32 & 28 &     \\
		&                                                 & balanced&     & 4.66 & 2.72 & 1.09 & 52 & 45 & 39 & 19\\
		& & medical &     &      & 3.73 &  &   60 & 59 & 58 &    \\  \midrule  
		\multirow{3}{*}{40}&  \multirow{3}{*}{5} & general          &6.89 & 3.95 & 2.83  & &46 & 38 & 32 &    \\
		& & balanced&     &      & 6.00 & 1.46 & 64 & 52 & 48 & 26 \\
		& & medical&     &      & 4.58 &   &   &    & 45 &      \\
		\bottomrule
	\end{tabular}
\end{table}

\paragraph{Percentage of split options utilized}
\cref{tab:decrease_in_oper_cost} shows the average percentages of utilized splits in columns eight to eleven. We observe that while task-splitting is beneficial in all considered scenarios, typically less than 50\% of the splittable visits are split for the instances considered. This confirms the relevance of optimizing task-splitting decisions instead of deciding which visits to split in advance. This allows to identify visits for which task-splitting is effective and to improve HHC planning without performing additional visits to a patient that are not beneficial. 
Comparing the different scenarios, the results indicate that an increase in the skill levels of the available staff (from practically trained staff to only medically trained staff) results in a decrease in the number of used splits because higher-skilled staff are capable of performing more tasks. Moreover, we observe a tendency for additional task-splitting when visits require more medical tasks (i.e., moving from general visit requirements to medical requirements) for a given staff composition (practically trained staff up to medically trained staff).

\smallskip

The main conclusions drawn from the results discussed in this section are that introducing task-splitting possibilities can reduce the operational costs for an HHC provider, help to deal with staff shortages, and help to ensure a better match between qualifications of caregivers and the tasks assigned to them. Our results also indicate that the full potential of task-splitting possibilities can only be used when task-splitting, routing, and scheduling decisions are combined (instead of using a sequential approach in which all splittable visits are split upfront). This can result in a larger number of utilized splits when visit requirements increase or staff qualifications decrease. 
Finally, our results also show that the benefits of task-splitting possibilities do not have to come at the price of overly dense schedules for caregivers, which could be detrimental to their work satisfaction. 

\smallskip 

The additional results provided in \ref{app:travel_time_minimization_res} indicate that the observed benefits of task-splitting do not depend on the specific choice of the objective function and carry over to the case of travel-time minimization, an objective that is frequently considered in the HHC literature.

\section{Conclusions}\label{sec:conclusion}
This is the first paper to formally examine the potential of incorporating task-splitting possibilities into home healthcare (HHC) planning. Home healthcare routing and scheduling problems (HHCRSPs) are widely studied due to the relevance of affordable and reliable HHC. Consequently, numerous problem variants have been proposed in the literature (e.g., \cite{OLADZADABBASABADY2023105829, BRAEKERS2016428}). Existing research focuses on the construction of efficient caregiver schedules, while typically accounting for varying qualifications and preferences among caregivers, specific service and time requests for patient visits, and occasionally temporal restrictions between multiple visits to the same patient. However, none of the papers studying home healthcare planning explores the potential of splitting patient visits consisting of multiple tasks into separate visits. 

This paper presents a novel problem variant of the HHCRSP, in which the decision whether certain patient visits should be split is integrated into the planning process. The objective is to design caregiver routes and schedules that minimize operational costs for home healthcare providers. The proposed model formulation allows split parts to have potentially wider visiting times, a shorter or longer combined duration, and to use lower-qualified caregivers than the unsplit visit.
Temporal dependencies can be imposed between all types of patient visits. Two novel integer linear programming formulations are presented (a MTZ formulation and a time-indexed variant) and an effective branch-and-bound solution framework is developed. The solution framework performance on the time-indexed variant is enhanced with pre-processing techniques and one of the two proposed heuristics, which is embedded as a primal and improvement heuristic in the branch-and-bound solution approach. The resulting solution method has potential to support decision making on task-splitting and scheduling in home healthcare. Considering the eventual effects on HHC schedules and routes, it selected only those task-splitting possibilities that contribute to an improved planning.

\smallskip

An extensive computational analysis examines the managerial and computational impacts of task-splitting. In the instances considered, optimized task-splitting allows to provide care with less caregivers or caregivers with less medical training and reduces operational costs for HHC providers. This enables a more optimal use of the available staff. Task-splitting can also help to ensure that caregivers spend relatively more time on care tasks commensurate with their qualification level. Another observation on our benchmark instances is that task-splitting can also have benefits when minimizing total travel time. 
From a computational point of view, the results demonstrate that introducing task-splitting poses a significant additional computational challenge due to the increase in routing options and synchronization requirements. A general-purpose solver was not able to find solutions for feasible instances within the time limit. A primal heuristic based on the MTZ formulation is capable of partially overcoming this issue; it found solutions for almost all of the considered test instances.

\smallskip 

Our research demonstrates the potential of the proposed formulations and solution framework and indicates the benefits of task-splitting. However, these findings could be strengthened by conducting an extensive pilot project in close collaboration with HHC providers. Such research could provide additional insights into the savings (in cost or working time) that can be achieved in practice, along with the conditions under which task-splitting can be embraced by staff and patients alike. In addition, one might discover the need to incorporate additional features into the model formulation, such as considering a multi-period scheduling horizon, patient preferences on limiting the number of different caregivers, patient preferences concerning the amount of time they need to remain at home due to (split) HHC visits, caregivers asking for workload balancing schedules, and multiple ways to split a visit (or vice versa, organize tasks into visits). 
Moreover, it may be useful to assess the impact of task-splitting for different instance characteristics not addressed in this work. These include, for example, visit duration, geographical area, the relation between the (qualification) characteristics of a splittable visit and both of its split parts, and the relation between the duration and required qualification type of a task/visit. 
Another avenue for further research is the development of solution methods capable of solving larger problem instances within a limited computation time. Although the methodology proposed in the current paper is capable of handling a wide variety of temporal restrictions associated with task-splitting, its computational performance for the test instances consisting of 40 visits is too high to support the daily planning of larger (home healthcare) service providers. Consequently, we see a potential for the development of alternative exact or heuristic solution frameworks that address synchronization constraints in a multitude of industries and problem settings. 
Since real-life travel and service times can be non-deterministic, and demand can be dynamic, there is also the potential for the development of solution approaches addressing these aspects.

\section*{Acknowledgments}
We thank the Samenwerkende Universitaire Rekenfaciliteiten (www.surf.nl) for allowing us to use the National Supercomputer Snellius. 

\section*{Data availability}
The set of benchmark instances is available at \url{https://github.com/LoekvMontfort/Task-splitting-in-home-healthcare-routing-and-scheduling-data}.

\bibliographystyle{chicago} 
\bibliography{bib/PaperHHCRSP-TS.bib}

\appendix

\section{Adjustments to temporal dependencies due to pre-processing}\label{app:adjust:temp:dep}

The pre-processing step of \cref{sec:prepro}, replacing constraints~\eqref{eq:ti:syn_12} and \eqref{eq:ti:syn_21} with \eqref{pre_proc:syn}, requires the addition of extra temporal dependencies if at least one visit has multiple temporal dependencies with other visits. 
Consider a visit $v \in V$ that has a temporal restriction with visit $u \in V$ and visit $w \in V$, i.e. $\{u,v\} \in D$ and $\{v,w\} \in D$. Let $t^{*}_{v}$ be the arrival time and $t^{\prime}_{v}$ be the departure time at visit $v$, where we require $t^{*}_{v} \neq t^{\prime}_{v}$. Before the pre-processing step of removing redundant travel arcs, time $t^{\prime}_{v}$ had to satisfy the temporal restrictions with $u$ and $w$. Now each time within interval $[t^{*}_{v}, t^{\prime}_{v}]$ can be used. As a result, a solution where temporal dependency $\{u,v\}$ requires $t^{*}_{v}$ while temporal dependency $\{v,w\}$ requires time $t^{\prime}_{v}$ is feasible. Since $v$ cannot start at two different moments in time, such solutions must be excluded. Specifically, we have to ensure that there is at least one time point within $[t^{*}_{v}, t^{\prime}_{v}]$ that satisfies the temporal dependence with $u$ and with $w$. This is done by introducing temporal dependencies between $u$ and $w$, specifically those that were previously implicitly ensured. 

\smallskip 

Before specifying the temporal restrictions between $u$ and $w$, we first recap the explicit temporal restriction between the starting times of $u$ and $v$ and $v$ and $w$. We denote the starting time of visit $u$, $v$, and $w$ with $\tu$, $\tv$, and $\tw$, respectively. Depending on the value of precedence variable $p_{uv}$, the following restrictions arise from temporal dependency $\{u,v\}$:

\begin{align}
    \mbox{if $p_{uv}=1$:}\quad & \tv \in [\tu+\delta^{\mathrm{min}}_{uv},\tu +\delta^{\mathrm{max}}_{uv}] \label{eq:app:uv:1:uv}\\  
    \Leftrightarrow\ & \tu \in [\tv-\delta^{\mathrm{max}}_{uv}, \tv-\delta^{\mathrm{min}}_{uv}] \\
    \mbox{if $p_{uv}=0$:}\quad &  \tu \in [\tv+\delta^{\mathrm{min}}_{vu},\tv +\delta^{\mathrm{max}}_{vu}] \label{eq:app:uv:0:vu}\\  \Leftrightarrow\ & \tv \in [\tu-\delta^{\mathrm{max}}_{vu}, \tu-\delta^{\mathrm{min}}_{vu}] \quad 
       \label{eq:app:uv:0:uv}
\end{align}

Similarly, depending on the value of $p_{vw}$, temporal dependency $\{v,w\}$ induces constraints  

\begin{align}
    \mbox{if $p_{vw}=1$:}\quad &   \tw \in [\tv+\delta^{\mathrm{min}}_{vw},\tv +\delta^{\mathrm{max}}_{vw}] \label{eq:app:vw:1:vw}\\  
  \Leftrightarrow\ & \tv \in [\tw-\delta^{\mathrm{max}}_{vw}, \tw-\delta^{\mathrm{min}}_{vw}]
       \label{eq:app:vw:1:wv}\\
    \mbox{if $p_{uv}=0$:}\quad &    \tv \in [\tw+\delta^{\mathrm{min}}_{wv}, \tw+\delta^{\mathrm{max}}_{wv}] \label{eq:app:vw:0:wv} \\ 
   \Leftrightarrow\ & \tw \in  [\tv-\delta^{\mathrm{max}}_{wv}, \tv -\delta^{\mathrm{min}}_{wv}] 
       \label{eq:app:vw:0:vw}
\end{align}

These temporal restrictions between $u$ and $v$ and $v$ and $w$ naturally induce temporal restrictions between $u$ and $w$, since $\tu$ and $\tw$ are linked with each other through $\tv$. Specifically, for each combination of a value of $p_{uv}$ and $p_{vw}$, the bounds on the starting times of $v$ with respect to $u$ and $w$ can be combined to derive the permitted difference between starting times $\tu$ and $\tw$ of visits $u$ and $w$, respectively. The various combinations of $p_{uv}$ and $p_{vw}$ result in the following restrictions:
\begin{itemize}
    \item Case $p_{uv}=1$ and $p_{vw}= 1$: 
    by using bounds \eqref{eq:app:uv:1:uv} on the starting time of $v$ with respect to $u$ in \eqref{eq:app:vw:1:vw}, we observe the minimum and maximum difference of the starting time $\tw$ with respect to $\tu$. Specifically, $\tw \geq \tw$ is required with $\tw \in [\tu+\delta^{\mathrm{min}}_{uv} + \delta^{\mathrm{min}}_{vw}, \tu+\delta^{\mathrm{max}}_{uv}+\delta^{\mathrm{max}}_{vw}]$.
    \item Case $p_{uv} = 1$ and $p_{vw}= 0$: 
    by using bounds \eqref{eq:app:uv:1:uv} on the starting time of $v$ with respect to $u$, but this time in \eqref{eq:app:vw:0:vw}, we observe the temporal relation $\tw \in [\tu+\delta^{\mathrm{min}}_{uv}-\delta^{\mathrm{max}}_{wv},  \tu+\delta^{\mathrm{max}}_{uv}-\delta^{\mathrm{min}}_{wv}]$. Then, splitting the results for the two starting orders of $u$ and $w$ gives:  
    \begin{itemize}
        \item $\tw \in [\tu+\max\{\delta^{\mathrm{min}}_{uv}-\delta^{\mathrm{max}}_{wv},0\},  \tu+\delta^{\mathrm{max}}_{uv}-\delta^{\mathrm{min}}_{wv}]$ for $\tw \geq \tu$, whose starting order is only possible when $\delta^{\mathrm{max}}_{uv}-\delta^{\mathrm{min}}_{wv} \geq 0$.
        \item $\tu \in [\tw+\max\{\delta^{\mathrm{min}}_{wv}-\delta^{\mathrm{max}}_{uv},0\}, \tw+\delta^{\mathrm{max}}_{wv}-\delta^{\mathrm{min}}_{uv}]$, whose starting order requires $\delta^{\mathrm{max}}_{wv}-\delta^{\mathrm{min}}_{uv}$.
    \end{itemize}
    \item Case $p_{uv} =0$ and $p_{vw} = 1$: using bounds \eqref{eq:app:uv:0:uv} on the starting time of visit $v$ with respect to $u$ in \eqref{eq:app:vw:1:vw}, we observe that $\tw \in [\tu - \delta^{\mathrm{max}}_{vu}+\delta^{\mathrm{min}}_{vw}, \tu- \delta^{\mathrm{min}}_{vu}+\delta^{\mathrm{max}}_{vw}]$. Splitting the results for the two starting orders of visits $u$ and $w$ once again, we arrive at:
    \begin{itemize}
        \item $\tw \in [\tu + \max\{- \delta^{\mathrm{max}}_{vu}+\delta^{\mathrm{min}}_{vw},0\}, \tu- \delta^{\mathrm{min}}_{vu}+\delta^{\mathrm{max}}_{vw}]$ for $\tw \geq \tu$, whose starting order is only feasible when $- \delta^{\mathrm{min}}_{vu}+\delta^{\mathrm{max}}_{vw} \geq 0$.
        \item $\tu \in [\tw + \max \{-\delta^{\mathrm{max}}_{vw}+\delta^{\mathrm{min}}_{vu},0\}, \tw-\delta^{\mathrm{min}}_{vw}+\delta^{\mathrm{max}}_{vu}]$, whose starting order can only exist when $-\delta^{\mathrm{min}}_{vw}+\delta^{\mathrm{max}}_{vu} \geq 0$.
    \end{itemize} 
    \item Case $p_{uv} = 0$ and $p_{vw}= 0$: use restrictions \eqref{eq:app:vw:0:wv} on the starting time of visit $v$ with respect to $w$ in \eqref{eq:app:uv:0:vu} to derive that $\tu \geq \tw$ and $\tu \in [\tw+\delta^{\mathrm{min}}_{wv} + \delta^{\mathrm{min}}_{vu}, \tw+\delta^{\mathrm{max}}_{wv}+\delta^{\mathrm{max}}_{vu}]$. 
\end{itemize}
\cref{tab:temp_dep_app} summarizes the induced temporal restrictions between $u$ and $w$ for each feasible combination of $p_{uv}$ and $p_{vw}$. Note that in case of an identified infeasible starting order of $u$ and $w$, we must replace the corresponding minimum and maximum difference of starting times by $\Tmax$ instead of the values denoted in the table. 

\smallskip  

As a result of the pre-processing, the above derived implicit induced temporal dependencies no longer hold. Consequently, solutions become feasible for which temporal dependencies $\{u,v\}$ and $\{w,v\}$ use non-overlapping intervals of $[t^{*}_{v}, t^{\prime}_{v}]$. To exclude such solutions, we reintroduce the implicit induced temporal dependencies between $u$ and $w$ explicitly. This ensures that there is a point in $[t^{*}_{v}, t^{\prime}_{v}]$ that satisfies both temporal restrictions $\{u,v\}$ and $\{v,w\}$, and thereby feasibility of the temporal restrictions between $u$, $v$, and $w$. 

To enforce the temporal restriction between $u$ and $w$ explicitly, we add an adjusted version of synchronization constraints \eqref{pre_proc:syn} for each feasible combination of $p_{uv}$ and $p_{vw}$. We use values $\delta^{\cdot}_{\cdot,p_{uv},p_{vw}}$ from \cref{tab:temp_dep_app}, and a possibly common precedence variable $p_{uw}$ in each constraint. Additionally, we add a penalty term on the right-hand side of each constraint to render the synchronization constraints redundant for other combinations of $p_{uv}$ and $p_{vw}$ (e.g., penalty $((1-p_{uv}) + p_{vw})$ for $p_{uv}=1$ and $p_{vw}=0$). 
For example, considering precedence values $p_{uv}=1$ and $p_{vw}=0$, the added constraints are 
\begin{subequations}
\begin{equation}
    p_{uw} + \sum_{q\in Q_u} X^q(\delta^{-}(u^{t}) \setminus \LAwait) + \sum_{\substack{a= (w^{\ell}, z^{m}) \in \LAcare : \\ \ell < t + \delta_{uw,1,0}^{\textrm{min}}}} \sum_{q\in Q_{w} \cap Q_{z}} X_a^{q} \leq  2 + ((1-p_{uv}) + p_{vw}) \label{pre_proc:syn:app:12a}
\end{equation}
for each $u^{t} \in \LNodesV{u} : |\delta^{-}(u^{t}) \setminus \LAwait| \geq 1$,  
\begin{equation}
    p_{uw} + \sum_{q\in Q_u} X^q(\delta^{+}(u^{t}) \cap \LAcare) + \sum_{\substack{a= (z^{m}, w^{\ell}) \in \LArcs \setminus \LAwait : \\ \ell > t + \delta_{uw,1,0}^{\textrm{max}}}} \sum_{q\in Q_{z} \cap Q_{w}} X_{a}^{q} \leq  2  + ((1-p_{uv}) + p_{vw}) \label{pre_proc:syn:app:12b}
\end{equation}
for all $u^{t} \in \LNodesV{u} : \ |\delta^{+}(u^{t}) \setminus \LAwait| \geq 1$, 
\begin{equation}
      \sum_{q\in Q_w} X^q(\delta^{-}(w^{t}) \setminus \LAwait)  + \sum_{\substack{a= (u^{\ell}, z^{m}) \in \LAcare : \\ \ell < t + \delta_{wu,1,0}^{\textrm{min}}}} \sum_{q\in Q_{u} \cap Q_{z}} X_a^{q}
    \leq 1 + p_{uw}    + ((1-p_{uv}) + p_{vw}) \label{pre_proc:syn:app:21a}\\
\end{equation}
for every $w^{t} \in \LNodesV{w} : |\delta^{-}(w^{t}) \setminus \LAwait| \geq 1$, and 
\begin{equation}
    \sum_{q\in Q_w} X^q(\delta^{+}(w^{t}) \cap \LAcare) + \sum_{\substack{a= (z^{m}, u^{\ell}) \in \LArcs \setminus \LAwait : \\ \ell > t + \delta_{wu,1,0}^{\textrm{max}}}} \sum_{q\in Q_{z} \cap Q_{u}} X_{a}^{q}  \leq 1 + p_{uw}   + ((1-p_{uv}) + p_{uv})  
\label{pre_proc:syn:app:21b}
\end{equation}
for each $w^{t} \in \LNodesV{w} : |\delta^{+}(w^{t}) \setminus \LAwait| \geq 1$.
\end{subequations}

\smallskip 

Note that if $u$ must start at the same time as $v$, it is sufficient to specify additional temporal dependencies for only one pair $\{u,w\}$. Specifically, we simply pick $\delta_{uw}^{\mathrm{min}}= \delta_{vw}^{\mathrm{min}}, \; \delta_{uw}^{\mathrm{max}}= \delta_{vw}^{\mathrm{max}}, \; \delta_{wu}^{\mathrm{min}}= \delta_{wv}^{\mathrm{min}}, \; \delta_{wu}^{\mathrm{max}}= \delta_{wv}^{\mathrm{max}}$, and do not require a penalty term. 

\begin{table}
\centering
\caption{Temporal difference between the starting times of two visits $u$ and $w$ if both have a temporal restriction with visit $v$.}
\label{tab:temp_dep_app}
\begin{tabular}{cccccc}
	\toprule
 & & \multicolumn{2}{c}{$\tw \geq \tu$} & \multicolumn{2}{c}{$\tu \geq \tw$} \\
\cmidrule(lr){3-4} \cmidrule(lr){5-6} 
$p_{uv}$ & $p_{vw}$ &$\delta^{\mathrm{min}}_{uw,p_{uv},p_{vw}}$  & $\delta^{\mathrm{max}}_{uw,p_{uv},p_{vw}}$  & $\delta^{\mathrm{min}}_{wu,p_{uv},p_{vw}}$  & $\delta^{\mathrm{max}}_{wu,p_{uv},p_{vw}}$    \\  
\midrule 
1 & 1 & $\delta^{\mathrm{min}}_{uv} + \delta^{\mathrm{min}}_{vw}$ & $\delta^{\mathrm{max}}_{uv}+\delta^{\mathrm{max}}_{vw}$ & $\Tmax$ & $\Tmax$ \\ 
1 & 0 & $\max\{\delta^{\mathrm{min}}_{uv}-\delta^{\mathrm{max}}_{wv},0\}$ & $\delta^{\mathrm{max}}_{uv}-\delta^{\mathrm{min}}_{wv}$ & $\max\{\delta^{\mathrm{min}}_{wv}-\delta^{\mathrm{max}}_{uv},0\}$  & $\delta^{\mathrm{max}}_{wv}-\delta^{\mathrm{min}}_{uv}$ \\
0 & 1 & $\max\{- \delta^{\mathrm{max}}_{vu}+\delta^{\mathrm{min}}_{vw},0\}$ & $- \delta^{\mathrm{min}}_{vu}+\delta^{\mathrm{max}}_{vw}$ & $\max \{-\delta^{\mathrm{max}}_{vw}+\delta^{\mathrm{min}}_{vu},0\}$ & $-\delta^{\mathrm{min}}_{vw}+\delta^{\mathrm{max}}_{vu}$ \\  
0 & 0 & $\Tmax$ & $\Tmax$ & $\delta^{\mathrm{min}}_{wv} + \delta^{\mathrm{min}}_{vu}$ & $\delta^{\mathrm{max}}_{wv}+\delta^{\mathrm{max}}_{vu}$\\
\bottomrule
\end{tabular}
\end{table}

\section{Additional computational results}
This section discusses additional computational results. \ref{app:perf_heuristics} elaborates on the performance of the different heuristic procedures, while \ref{app:devision_working_time} focuses on the division of working time among the different caregiver types with and without task-splitting. 
\subsection{Performance of heuristics}\label{app:perf_heuristics} 
\begin{figure}[h]
\centering
\includegraphics[width=1.0\linewidth]{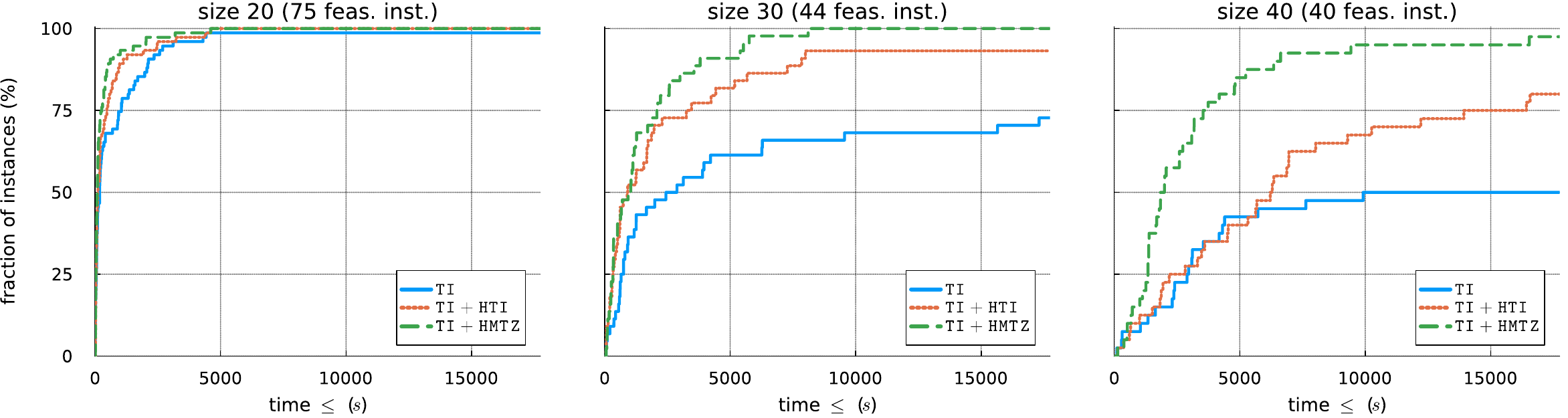}
\caption{Runtimes [s] until the first solution for \TI, \TIHTI, and \TIHMTZ. 
}
\label{fig:time_first_sol_3_approaches}
\end{figure}

\begin{figure}[h]
\centering
\includegraphics[width=1.0\linewidth]{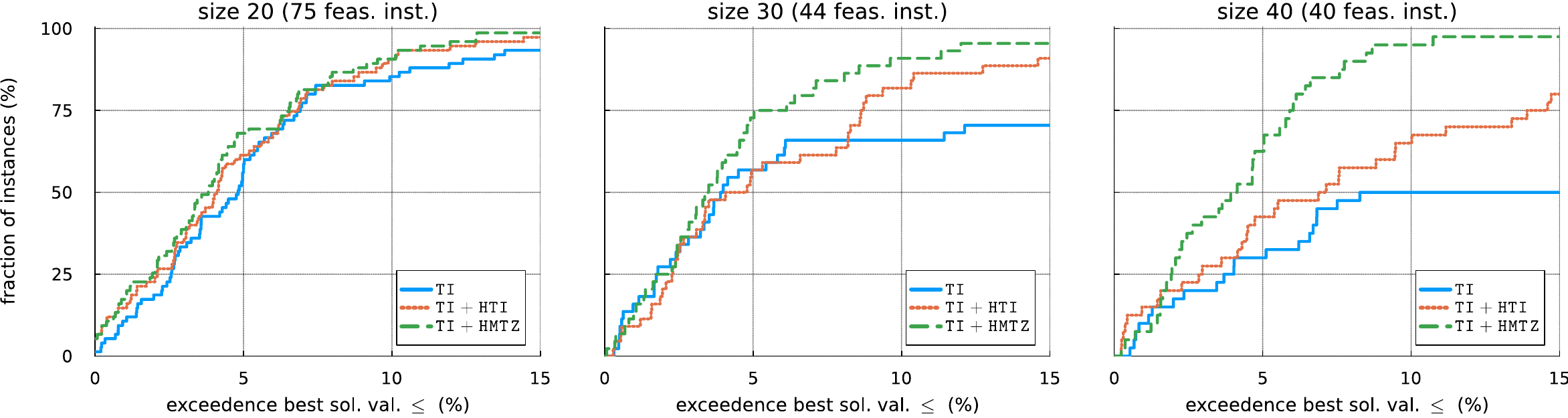}
\caption{Relative gaps between the first solution found by \TI, \TIHTI, and \TIHMTZ and overall best known solution. 
}
\label{fig:perc_first_sol_above_best_sol}
\end{figure}

\cref{fig:time_first_sol_3_approaches} and \cref{fig:perc_first_sol_above_best_sol} show the relative amounts of instances for which the first solution is found within a given time (\cref{fig:time_first_sol_3_approaches}) 
and for which the relative difference between the objective value of this first solution and the overall best-known solution for that instance is below a certain threshold (\cref{fig:perc_first_sol_above_best_sol}).
The largest value on the x-axis of each plot in \cref{fig:time_first_sol_3_approaches} corresponds to the maximum runtime of 18\,000 seconds. If the plot for a particular method does not reach 100\%, then this signals that no feasible solution was found for the remaining fraction of the instances. In \cref{fig:perc_first_sol_above_best_sol} the results are truncated at 15\%, which means that a relative gap exceeded 15\% or no feasible solution was found for the remaining instances of a particular method if the fraction of instances does not reach 100\%.

The results of \cref{fig:time_first_sol_3_approaches} show that solutions are found earlier when using the heuristic procedures and that the initial solutions they find are clearly better than the initial solutions of variant \TI obtained by Gurobi's general-purpose heuristics. These positive effects are more pronounced for \TIHMTZ than for \TIHTI and increase with an increasing instance size. 
An in-depth analysis of the results showed that the first feasible solution is found by a primal heuristic in 47.5\% and 87.5\% of all instances of size 40 in the case of \TIHTI and \TIHMTZ, respectively. Variant \TIHTI with the time-indexed primal heuristic required an average of 455 attempts, each taking 0.7 seconds, until a feasible solution was found from a fractional solution. In contrast, variant \TIHMTZ based on the MTZ formulation only required 61 calls lasting on average 0.4 seconds. 
Given the clear benefits of \TIHMTZ over the other two alternatives, we will concentrate on this variant in the following paragraph. 

\smallskip

With regard to the quality of initial solutions, we observe that 75\% and 67.5\% of those found by \TIHMTZ are, at most, 5\% worse than the best-known final solutions (among the \TI, \TIHTI, and \TIHMTZ methods) for instances of size 30 and 40, respectively. 
The improvement step contributes significantly to their quality. Focusing on the instances of size 40, an initial solution found by the MTZ-based primal heuristic could be enhanced by the improvement heuristic in 94.3\% of the cases, reducing the objective value by 6.0\% on average. The time required to apply the improvement heuristic on instances with size 40 is usually below 30 seconds (on average 94 seconds due to one outlier reaching the time limit). We note that for instances of size 40, the improvement heuristic also consistently enhanced initial solutions found by Gurobi's general-purpose heuristics, in which case its average objective value reduction was equal to 0.9\%. Furthermore, on all instances of size 40, it enhanced 56.9\% of the solutions found by Gurobi after the initial one. The average relative cost decrease was 0.2\% and the corresponding average runtime 0.3 seconds.

\subsection{Division of working time over the different caregiver types}\label{app:devision_working_time}
Tables \ref{tab: work_day_level1}, \ref{tab: work_day_level2}, and \ref{tab: work_day_level3} provide an overview of the average allocation of total working time among the three different types of caregivers for the instances considered, with and without task-splitting. Only instances feasible without task-splitting are considered in this analysis. In every scenario except for one, the percentage of working time assigned to level 3 caregivers decreases because of the introduction of task-splitting. Conversely, we observe that the percentage of total working time assigned to level 1 caregiver increases in each scenario. 

\begin{table}[H]
\centering
\small 
\caption{Average percentage of the total working time assigned to level 1 caregivers.}
\begin{tabular}{rrrcccccccc}
\toprule 
& & & \multicolumn{8}{c}{average percentage of working time assigned to level 1 caregivers (\%)} \\
\cmidrule(lr){4-11}
& & & \multicolumn{2}{c}{prac. train. staff} & \multicolumn{2}{c}{moder. train. staff} & \multicolumn{2}{c}{medic. train. staff} & \multicolumn{2}{c}{only medic. staff} \\
\cmidrule(lr){4-5} \cmidrule(lr){6-7} \cmidrule(lr){8-9} \cmidrule(lr){10-11}
size & $\#$ & visit req. & no split & split & no split & split & no split & split & no split & split \\ 
	\midrule
\multirow{3}{*}{20} & \multirow{3}{*}{10} & general &  46.6 & 53.3 & 34.5 & 37.4 & 32.6 & 34.5 &   &    \\
                                                & & balanced&  41.0 & 42.5 & 33.3 & 35.7 & 29.5 & 33.2 & 0.0 & 0.0   \\
                                                & &  medical&       &      &      &      & 19.5 & 32.0 &   &   \\  \midrule
\multirow{3}{*}{30} & \multirow{3}{*}{5}  & general &  50.2 & 57.1 & 22.3 & 24.2 & 23.5 & 24.2 &   &   \\ 
&                                                 & balanced&       &      & 21.1 & 23.4 & 20.8 & 22.8 & 0.0 & 0.0  \\
                                                & & medical &       &      &      &      & 19.7 & 22.6 &   &    \\ \midrule 
\multirow{3}{*}{40}&  \multirow{3}{*}{5} & general &  50.0 & 57.9 & 33.3 & 35.2 & 33.5 & 35.1 &   &   \\
                                                & & balanced&       &      &      &      & 28.8 & 33.6 & 0.0 & 0.0  \\
                                                 & & medical&       &      &      &      & 25.6 & 31.0 &   &    \\
                                                 \bottomrule
\end{tabular} 
\label{tab: work_day_level1}
\footnotesize
\end{table} 
\begin{table}[H]
\centering
\small 
\caption{Average percentage of the total working time assigned to level 2 caregivers.}
\begin{tabular}{rrrcccccccc}
\toprule 
& & & \multicolumn{8}{c}{average percentage of working time assigned to level 2 caregivers (\%)} \\
\cmidrule(lr){4-11}
& & & \multicolumn{2}{c}{prac. train. staff} & \multicolumn{2}{c}{moder. train. staff} & \multicolumn{2}{c}{medic. train. staff} & \multicolumn{2}{c}{only medic. staff} \\
\cmidrule(lr){4-5} \cmidrule(lr){6-7} \cmidrule(lr){8-9} \cmidrule(lr){10-11}
size & $\#$ & visit req. & no split & split & no split & split & no split & split & no split & split \\ 
	\midrule
\multirow{3}{*}{20} & \multirow{3}{*}{10} & general &24.9 & 18.8 & 36.4 & 33.7 & 27.9 & 30.2 &   &     \\
                                                & & balanced&32.0 & 29.2 & 32.6 & 30.3 & 29.4 & 30.3 & 0.0 & 0.0  \\
                                                & &  medical&     &      &      &      & 26.2 & 20.2 &   &  \\  \midrule  
\multirow{3}{*}{30} & \multirow{3}{*}{5}  & general &20.7 & 19.9 & 50.4 & 52.6 & 39.7 & 42.2 &   &   \\ 
&                                                 & balanced&     &      & 44.5 & 48.4 & 39.1 & 40.1 & 0.0 & 0.0\\
                                                & & medical &     &      &      &      & 31.2 & 35.6 &   & \\ \midrule 
\multirow{3}{*}{40}&  \multirow{3}{*}{5} & general &22.8 & 20.0 & 39.6 & 42.6 & 31.9 & 33.3 &   &\\
                                                & & balanced&     &      &      &      & 30.2 & 31.7 & 0.0 & 0.0  \\
                                                 & & medical&     &      &      &      & 23.7 & 24.9 &   &   \\ 
                                                 \bottomrule
\end{tabular} 
\label{tab: work_day_level2}
\footnotesize
\end{table} 
\begin{table}[H]
\centering
\small 
\caption{Average percentage of the total working time assigned to level 3 caregivers.}
\begin{tabular}{rrrcccccccc}
\toprule 
& & & \multicolumn{8}{c}{average percentage of working time assigned to level 3 caregivers (\%)} \\
\cmidrule(lr){4-11}
& & & \multicolumn{2}{c}{prac. train. staff} & \multicolumn{2}{c}{moder. train. staff} & \multicolumn{2}{c}{medic. train. staff} & \multicolumn{2}{c}{only medic. staff} \\
\cmidrule(lr){4-5} \cmidrule(lr){6-7} \cmidrule(lr){8-9} \cmidrule(lr){10-11}
size & $\#$ & visit req. & no split & split & no split & split & no split & split & no split & split \\ 
	\midrule
\multirow{3}{*}{20} & \multirow{3}{*}{10} & general  &28.4 & 28.0 & 29.1 & 28.8 & 39.5 & 35.3 &     &        \\
                                                & & balanced&27.0 & 28.3 & 34.1 & 34.0 & 41.1 & 36.5 & 100.0 & 100.0  \\
                                                & &  medical&     &      &      &      & 54.4 & 47.8 &     &     \\  \midrule 
\multirow{3}{*}{30} & \multirow{3}{*}{5}  & general &29.1 & 23.0 & 27.3 & 23.2 & 36.8 & 33.6 &     &  \\ 
&                                                 & balanced&     &      & 34.4 & 28.2 & 40.1 & 37.1 & 100.0 & 100.0\\
                                                & & medical &     &      &      &      & 49.1 & 41.8 &     &\\ \midrule  
\multirow{3}{*}{40}&  \multirow{3}{*}{5} & general  &27.1 & 22.0 & 27.0 & 22.2 & 34.6 & 31.7 &     & \\
                                                & & balanced&     &      &      &      & 41.0 & 34.8 & 100.0 & 100.0  \\
                                                 & & medical&     &      &      &      & 50.7 & 44.1 &     & \\ 
                                                 \bottomrule
\end{tabular} 
\label{tab: work_day_level3}
\footnotesize
\end{table} 

\section{Comparing the impact of task-splitting for operational costs minimization and travel time minimization}\label{app:travel_time_minimization}

This section discusses whether the potential benefits of task-splitting depend on the selected objective of operational costs minimization. It compares computational results for operational cost minimization with those for travel-time minimization to ascertain whether task-splitting can be relevant for both objectives. \ref{app:travel_time_minimization_prob_var} introduces the problem variant focused on travel time minimization together and the required modification for the mathematical formulations, while \ref{app:travel_time_minimization_res} compares the computational results for the two objective functions.

\subsection{Problem variant focusing on travel time minimization}\label{app:travel_time_minimization_prob_var}
The \probabbr includes the wage of the entire working shift of a caregiver in the objective function. In practice, travel time between patients is typically defined as working time (see, for example, the Dutch labor agreement; \citetalias{COAVVT2023}). 
To reflect the operational costs of a caregiver for an HHC provider, the \probabbr incorporates a caregiver's wage from the start of their first visit until the completion of their last visit. 
However, related research focuses, for example, on the total travel time of caregivers (see, e.g., \cite{PAHLEVANI2022102878, Frifita2020}). 
To explore the impact of task-splitting for different objective functions, we also consider a variant of the \probabbr, denoted as \probabbrtt, whose objective is to minimize the total travel time of all caregivers. The \probabbrtt neglects visit time, waiting time, and wage differences that are included in the \probabbr and might effect the impact of task-splitting.

\medskip 

Both the MTZ formulation~\eqref{eq:mtz} and the time-indexed formulation~\eqref{eq:ti} have to be adjusted for this new problem variant. 
A MTZ formulation for the \probabbrtt is derived from formulation~\eqref{eq:mtz} by replacing the objective function~\eqref{eq:mtz:obj} by $\sum_{a=(u,v) \in A}\sum_{q \in Q_{u} \cap Q_{v}} t_{uv}x_{a}^{q}$, which sums up the travel time of all travels performed. In addition, variables $k^q_v$ and $f^q_v$ are removed $\forall v\in V$, $\forall q\in Q_v$, as well as constraints~\eqref{eq:mtz:starttimes}--\eqref{eq:mtz:finishtimes:forcezero}.
Similarly, a time-indexed formulation for the \probabbrtt is obtained from \eqref{eq:ti} by inserting the objective $\sum_{a=(u^{t},v^\ell) \in \Acare} \sum_{q \in Q_{u} \cap Q_{v}} t_{uv} X_{a}^{q}$ in place of \eqref{eq:ti:obj}, which represents the total travel time of the caregivers.

\subsection{Computational results}\label{app:travel_time_minimization_res}
In this section we will assess to what extent the observed benefits of task-splitting depend on the specific choice of the objective function. To this end, we compare the results obtained using \TIHMTZ for operational cost minimization (\probabbr) with those for travel time minimization (\probabbrtt). We will first focus on reducing objective function values, then analyze the effect of task-splitting on staff schedules, and lastly compare the percentage of split options utilized in the solutions.

\paragraph{Objective function value}
\cref{tab:decrease_in_objective_value_both_app} presents the average reduction in objective function values thanks to task-splitting for operational cost minimization (discussed in \cref{sec:results:impact}) and travel time minimization. The objective of the \probabbrtt does not consider wage differences between the different caregiver types. Moreover, additional travel time due to task-splitting cannot be compensated with reduced care or shorter waiting times. \cref{tab:decrease_in_objective_value_both_app} illustrates that task-splitting can be interesting for HHC providers when minimizing travel time. Splitting visits, for example, can create additional routes that contribute to a solution with a lower total travel time. 
Although the average reduction of the objective function value is slightly smaller compared to operational cost minimization (the weighted averages of the results reported in \cref{tab:decrease_in_objective_value_both_app} are equal to 3.7\% and 4.0\%), we consistently observe gains in all considered scenarios that are sometimes even larger than those for the \probabbr. The results for the benchmark instances clearly indicate that task-splitting can not only be relevant when minimizing the operational costs and considering wage differences between different types of caregivers. It can also be relevant when minimizing travel time. 

\begin{table}
\centering
\caption{Average objective value decrease due to task-splitting when minimizing operational costs (o.c.) or travel time (t.t.).}
\label{tab:decrease_in_objective_value_both_app}
\small
\begin{tabular}{rrrcccccccc}
	\toprule
& & & \multicolumn{8}{c}{average decrease in objective value due to task-splitting (\%)} \\
\cmidrule(lr){4-11}
& & & \multicolumn{2}{c}{prac. train. staff} & \multicolumn{2}{c}{moder. train. staff} & \multicolumn{2}{c}{medic. train. staff} & \multicolumn{2}{c}{only medic. staff} \\
\cmidrule(lr){4-5} \cmidrule(lr){6-7} \cmidrule(lr){8-9} \cmidrule(lr){10-11}
size & $\#$ & visit req. & o.c. & t.t. & o.c. & t.t. & o.c. & t.t. & o.c. & t.t.  \\ 
	\midrule
\multirow{3}{*}{20} & \multirow{3}{*}{10} & general &3.83 & 4.80 & 2.39 & 0.45 & 4.33 & 3.39 &      &   \\
                                                & & balanced&4.68 & 0.39 & 4.26 & 1.93 & 5.97 & 1.45 & 1.16 & 1.50   \\
                                                & &  medical&     &      &      &      & 8.08 & 9.08 &      & \\ \midrule 
\multirow{3}{*}{30} & \multirow{3}{*}{5}  & general &7.26 & 6.72 & 3.37 & 2.60 & 2.23 & 1.40 &      &      \\
&                                                 & balanced&     &      & 4.66 & 9.88 & 2.72 & 2.22 & 1.09 & 0.64\\
                                                & & medical &     &      &      &      & 3.73 & 6.44 &      & \\  \midrule  
\multirow{3}{*}{40}&  \multirow{3}{*}{5} & general &6.89 & 6.22 & 3.95 & 7.54 & 2.83 & 1.77 &      &     \\
                                                & & balanced&     &      &      &      & 6.00 & 2.19 & 1.46 & 1.17 \\
                                                 & & medical&     &      &      &      & 4.58 & 8.55 &      &    \\
	\bottomrule
\end{tabular}
\footnotesize
\end{table}
\paragraph{Staff schedules}
To further explore the differences in solutions obtained for the two objectives, \cref{tab:perc_patient_time_both_objectives_app} reports the average percentage of time caregivers spend on patient care. To facilitate a fair comparison, the starting times of patient visits for solutions found in variant \probabbrtt are improved during post-processing to avoid unnecessary waiting time in the schedules.
We note that no feasible solution could be found for a single instance of size 40 when minimizing travel time. Consequently, this instance is excluded from this comparison.
The largest difference in the percentage of working time spent on patient care occurs if only medically trained staff is available. In these generated scenarios, minimizing travel time instead of operational costs results in caregiver schedules for which the percentage of working time spent on patient care is, on average, 8.5 percentage points lower. 
A summary of all considered benchmark instances indicates that the percentage of a workday that caregivers spend on providing care is, on average, four percentage points higher when minimizing operational costs. 
Since waiting time is not considered in the objective of the \probabbrtt, the length of the resulting schedules is not optimized. As a result, caregivers spend, on average, approximately 5.6\% of their workday waiting for the benchmark instances. In contrast, the schedules obtained with \probabbr, which minimizes operational costs, are much more compact and efficient for caregivers, requiring less working time for the same amount of care. 
In addition to potentially shorter workdays, which are interest for caregivers, the operational costs objective is also relevant for HHC providers. It can create time for additional care tasks and may also reduce (variable) labor costs.
Therefore, an operational cost objective seems to be a more natural option for both stakeholders. 

We note that the division of working time among the different caregiver types is different for both objective functions on the instances considered. When minimizing travel time instead of operational costs, a larger portion of the total working time is assigned to medically trained caregivers in a situation without task-splitting.
After task-splitting, the division of working time between the different caregiver types only remains roughly the same for travel time minimization.

\begin{table}
\centering
\caption{Average percentage of the total working time spent on patient care when minimizing the operational costs (o.c.) or travel time (t.t.). 
}
\label{tab:perc_patient_time_both_objectives_app}
\small
\begin{tabular}{rrrcccccccc}
	\toprule
& & & \multicolumn{8}{c}{average percentage of working time spent on patient care (\%)} \\ 
\cmidrule(lr){4-11}
& & & \multicolumn{2}{c}{prac. train. staff} & \multicolumn{2}{c}{moder. train. staff} & \multicolumn{2}{c}{medic. train. staff} & \multicolumn{2}{c}{only medic. staff} \\
\cmidrule(lr){4-5} \cmidrule(lr){6-7} \cmidrule(lr){8-9} \cmidrule(lr){10-11}
size & $\#$ & visit req. & o.c. & t.t. & o.c. & t.t. & o.c. & t.t. & o.c. & t.t.  \\ 
	\midrule
\multirow{3}{*}{20} & \multirow{3}{*}{10} & general &79.4 & 76.3 & 81.9 & 77.2 & 82.5 & 76.4 &      &       \\
                                                & & balanced&79.0 & 75.8 & 80.3 & 76.6 & 81.8 & 76.8 & 85.8 & 76.4   \\
                                                & &  medical&     &      &      &      & 80.1 & 76.4 &      &         \\  \midrule 
\multirow{3}{*}{30} & \multirow{3}{*}{5}  & general &83.9 & 83.5 & 85.8 & 82.9 & 86.5 & 82.4 &      &         \\
&                                                 & balanced&81.5 & 82.4 & 83.6 & 80.3 & 85.2 & 82.4 & 89.1 & 83.6   \\
                                                & & medical &76.6 & 72.8 & 80.8 & 76.2 & 83.8 & 82.6 &      &        \\  \midrule  
\multirow{3}{*}{40}&  \multirow{3}{*}{5} & general &84.0 & 81.8 & 86.4 & 82.0 & 87.9 & 83.1 &      &        \\
                                                & & balanced&84.2 & 83.2 & 83.7 & 81.5 & 86.1 & 82.8 & 90.4 & 80.8      \\
                                                 & & medical&     &      &      &      & 84.2 & 82.9 &      &    \\ 
	\bottomrule
\end{tabular}
\footnotesize
\end{table}

\paragraph{Percentage of split options utilized}

\cref{tab:perc_splits_both_objectives_app} shows the average percentages of utilized splits for the two considered problem variants. Only instances for which a feasible solution was found for both problem variants are considered in this analysis. 
A comparison of the results for the two objective functions reveals that a greater number of splits occur for the operational cost minimization objective. This may be attributed to the fact that, in the case of operational cost minimization, extra travel time can be compensated by a reduction in patient time and/or waiting time, thereby making task-splitting more beneficial. Wage differentiation also stimulates task-splitting, resulting in the execution of tasks by the cheapest, appropriately qualified caregiver. 
In turn, the comparably low number of splits for the travel time minimization objective suggests that even a few splits (see \cref{tab:perc_splits_both_objectives_app}) can have benefits for caregiver planning in terms of the number of instances that are feasible and the reduction in terms of objective value that can be achieved (see \cref{tab:number_of_feasible_instances} and \cref{tab:decrease_in_objective_value_both_app}).

\begin{table}
\centering
\caption{Average percentage of utilized splits when minimizing operational costs (o.c.) or travel time (t.t.).
}
\label{tab:perc_splits_both_objectives_app}
\small
\begin{tabular}{rrrcccccccc}
	\toprule
& & & \multicolumn{8}{c}{average percentage of utilized splits (\%)} \\ 
\cmidrule(lr){4-11}
& & & \multicolumn{2}{c}{prac. train. staff} & \multicolumn{2}{c}{moder. train. staff} & \multicolumn{2}{c}{medic. train. staff} & \multicolumn{2}{c}{only medic. staff} \\
\cmidrule(lr){4-5} \cmidrule(lr){6-7} \cmidrule(lr){8-9} \cmidrule(lr){10-11}
size & $\#$ & visit req. & o.c. & t.t. & o.c. & t.t. & o.c. & t.t. & o.c. & t.t.  \\ 
	\midrule
\multirow{3}{*}{20} & \multirow{3}{*}{10} & general &40 & 21 & 29 & 16 & 29 & 8  &    &     \\
                                                & & balanced&43 & 27 & 35 & 19 & 32 & 6  & 21 & 3   \\
                                                & &  medical&   &    &    &    & 47 & 18 &    &   \\  \midrule 
\multirow{3}{*}{30} & \multirow{3}{*}{5}  & general &44 & 13 & 32 & 6  & 28 & 6  &    & \\       
&                                                 & balanced&52 & 33 & 45 & 17 & 39 & 5  & 19 & 7     \\ 
                                                & & medical &60 & 47 & 59 & 37 & 58 & 11 &    &  \\ \midrule  
\multirow{3}{*}{40}&  \multirow{3}{*}{5} & general  &46 & 9  & 38 & 11 & 32 & 4  &    &         \\
                                                & & balanced&60 & 35 & 52 & 22 & 48 & 3  & 26 & 5    \\
                                                 & & medical&   &    &    &    & 45 & 15 &    &    \\ 
	\bottomrule
\end{tabular}
\footnotesize
\end{table}

\end{document}